\documentclass[12pt,reqno]{article}
mm \textheight=8.9in

\usepackage{amssymb}
\usepackage{amsmath}
\usepackage{amsthm}

\usepackage{color}
\usepackage{pb-diagram}
\newcommand{\tw}[3]{{$#1$}${\,\scriptscriptstyle {#2}}\atop\raise9pt\hbox{$\scriptstyle\tp$} ${$#3$}}
\newcommand{\st}[1]{\mbox{${\,\scriptscriptstyle {#1}}\atop\raise5.5pt\hbox{$*$}$}}
\newcommand{\btr}{\raise1.2pt\hbox{$\scriptstyle\blacktriangleright$}\hspace{2pt}}
\newcommand{\id}{\mathrm{id}}

\newcommand{\tp}{\otimes}

\newcommand{\Lc}{\mathcal{L}}
\newcommand{\A}{\mathcal{A}}

\newcommand{\Rg}{\mathcal{L}}
\newcommand{\Bg}{\mathfrak{B}}
\newcommand{\B}{\mathcal{B}}
\newcommand{\D}{\mathfrak{D}}

\newcommand{\E}{\mathfrak{E}}

\newcommand{\Cc}{\mathcal{C}}
\newcommand{\Mc}{\mathcal{M}}

\newcommand{\C}{\mathbb{C}}

\newcommand{\Ha}{\mathcal{H}}
\newcommand{\Ru}{\mathcal{R}}

\newcommand{\op}{{\bullet}}

\newcommand{\sst}{\scriptscriptstyle}
\newcommand{\U}{\mathcal{U}}
\newcommand{\F}{\mathcal{F}}
\newcommand{\ve}{\varepsilon}
\newcommand{\gm}{\gamma}

\newcommand{\la}{\lambda}

\newcommand{\End}{\mathrm{End}}

\newcommand{\Hom}{\mathrm{Hom}}
\newcommand{\Ob}{\mathrm{Ob}}

\newcommand{\tr}{\triangleright}

\newcommand{\btl}{\mbox{\raise1.1pt\hbox{$\scriptstyle\blacktriangleleft$}}}

\newcommand{\g}{\mathfrak{g}}

\newcommand{\h}{\mathfrak{h}}

\newcommand{\nn}{\nonumber}

\newcommand{\si}{\sigma}
\newcommand{\al}{\alpha}
\newcommand{\bt}{\beta}
\newcommand{\Mod}{\mathrm{Mod }\:}
\newcommand{\Modd}{\mathrm{Mod_0}\:}

\newcommand{\be}{\begin{eqnarray}}
\newcommand{\ee}{\end{eqnarray}}

\newtheorem{thm}{Theorem}[section]
\newtheorem{propn}[thm]{Proposition}
\newtheorem{lemma}[thm]{Lemma}
\newtheorem{corollary}[thm]{Corollary}

\theoremstyle{definition}
\newtheorem{remark}[thm]{Remark}
\newtheorem{definition}[thm]{Definition}
\newtheorem{example}[thm]{Example}
\newtheorem{remarks}[thm]{Remarks}
\newtheorem{examples}[thm]{Examples}

\newcount\prg

\newcommand{\parag}{\advance\prg by1 {\noindent\bf\thesection.\the\prg\hspace{6pt}}}
\newcommand{\select}[1]{\textcolor{red}{\bf\em #1}}

\begin{document}
\title{Quantum groupoids and dynamical categories\footnote{
This research is supported in part
by the Israel Academy of Sciences grant no. 8007/99-03,
the Emmy Noether Research Institute for Mathematics,
the Minerva Foundation of Germany,  the Excellency Center "Group
Theoretic Methods in the study of Algebraic Varieties"  of the Israel
Science foundation, and by the RFBR grant no. 03-01-00593.}}
\author{J. Donin$^\dagger$
\hspace{3pt} and A. Mudrov$^\ddag$}
\date{}
\maketitle
\begin{center}
$\dagger$,$\ddag${Department of Mathematics, Bar Ilan University, 52900 Ramat Gan,
Israel,\\
$\ddag$
Max-Planck Institut f$\ddot{\rm u}$r Mathematik, Vivatsgasse 7, D-53111 Bonn, Germany.}
\end{center}
\begin{abstract}
In this paper we realize the dynamical categories introduced in our previous paper
as categories of modules over bialgebroids; we study the
bialgebroids arising in this way.
We define quasitriangular structure on bialgebroids and
present examples of quasitriangular bialgebroids related to the dynamical categories.
We show that dynamical twists over an arbitrary  base
give rise to bialgebroid twists.

We prove  that the classical dynamical r-matrices
over an arbitrary base manifold are in one-to-one correspondence with
a special class of coboundary Lie bialgebroids.
\end{abstract}
{\small \underline{Key words}: Dynamical category, dynamical twist,
dynamical Yang-Baxter equation, bialgebroid, quantum groupoid, Lie bialgebroid.\\
\underline{AMS classification codes}: 17B37, 81R50.}
\maketitle
\tableofcontents
\section{Introduction}
In our recent paper \cite{DM1}, we introduced
a procedure of dynamization of monoidal categories.
The categorical approach
naturally led to a definition of the dynamical Yang-Baxter equation (DYBE), both classical and quantum,
over an arbitrary base\footnote{For an introduction to the theory of dynamical Yang-Baxter equation and the bibliography
see \cite{ESch1}}.
In this way, the constructions of twists from \cite{Xu2}
and \cite{EE1} acquired a categorical meaning.
The dynamical categories of  \cite{DM1} generalize the dynamical categories which were introduced
by Etingof and Varchenko \cite{EV2} for commutative cocommutative Hopf algebras.

In the framework of the categorical approach, we developed a fusion procedure which led to a
construction of dynamical twists.
Those dynamical twists were used for equivariant star product quantization
of vector bundles on the coadjoint orbits of reductive Lie groups, including
the algebras of functions (see also \cite{AL} and \cite{KMST}).
In a recent paper of Etingof and Enriquez \cite{EE2}, this fusion procedure was extended
further, including a class of infinite dimensional Lie algebras.

The goal of the present paper is to realize the dynamical categories as representations
of certain $\Lc$-bialgebroids, were $\Lc$ is a base algebra over a Hopf algebra $\Ha$
in the sense of Definition \ref{defBA}.
The simplest bialgebroid of this kind,
namely the smash product $\Lc\rtimes \Ha$, was introduced
in \cite{Lu}. It is interesting to note that bialgebroids of \cite{Lu} were considered over exactly the
same class of base algebras that was used for the definition of dynamical categories in \cite{DM1}.
In the present paper we link the theory of dynamical Yang-Baxter equations over a non-abelian base
with the bialgebroids of \cite{Lu}. We pursue a further study of those bialgebroids and their descenders.
In particular, we show that their certain quotients have a quasitriangular structure.

The infinitesimal analogs of bialgebroids are Lie bialgebroids.
We consider Lie bialgebroids which are quasi-classical limits of the bialgebroids
related to the dynamical categories. In this way we come to the most general definition
of dynamical r-matrix over an arbitrary base manifold as the space of dynamical parameters.
We show that the classical dynamical r-matrices are in one-to-one correspondence with
a special class of coboundary Lie bialgebroids.

In the present paper we obtain the following results.

We define quasitriangular structure and the
notion of universal R-matrix on bialgebroids.

We study quasitriangular bialgebroids (quantum groupoids) related to
the dynamical categories.

We give an interpretation to the antipode of \cite{Lu} as a an isomorphism between
two different bialgebroids over different bases.

We prove that a dynamical twist over an arbitrary base gives rise to a twist
of bialgebroids. This is a generalization of the results of \cite{Xu1}.

We present an example of a "dual" bialgebroid over an non-abelian base.

We define a classical dynamical r-matrix over a Poisson base algebra $\Lc_0$ as
a coboundary Lie bialgebroid  of a special type over $\Lc_0$.

The paper is organized as follows.

Section \ref{secDC} recalls the construction of dynamical categories over
an $\Ha$-base algebra $\Lc$ for some Hopf algebra, $\Ha$.

Section \ref{secHA}
contains basic definitions from the theory of bialgebroids.

Section \ref{secBH_R} introduces a bialgebroid extension of a quasitriangular
Hopf algebra  $\Ha$ by its quasi-commutative module algebra $\Rg$.
Therein we show that a certain quotient bialgebroid $\Ha_\Rg$ has a quasitriangular structure.

In Section \ref{secAntipode} we give an interpretation of Lu's antipode on
the smash product bialgebroid as an isomorphism between a pair of bialgebroids.
We prove that the antipode is carried over to the quotient quantum groupoid $\Ha_\Rg$.

Section \ref{secBTDC} establishes a relation between dynamical cocycles and bialgebroid twists.
We start from the trivial extension of the bialgebroid $\D\Ha_\Lc$, where $\D\Ha$ is the double
of $\Ha$, by a Hopf algebra $\U$ containing $\Ha$.
We show that the element $\Psi=\F\Theta$ built out of a dynamical cocycle $\F\in \U\tp \U\tp \Lc$
and a universal R-matrix $\Theta$ of the double,
is a twist of the bialgebroid $\U\tp \D\Ha_\Lc$.

Section \ref{secFDCR} realizes dynamical categories as representations of bialgebroids.

In Section \ref{secDQG} we present a "bimodule" algebra
over the tensor product bialgebroid $\U\tp \D\Ha_\Lc$ twisted by a dynamical twist.
It is, in fact, a bialgebroid and may be considered as a dynamical FRT algebra,
in case $\U$ is quasitriangular.

In Section \ref{secCDRM} we give the most general, to our knowledge, definition of the
classical dynamical r-matrix over arbitrary base. We prove that a classical dynamical r-matrix in the
sense of that definition is the same as a special coboundary Lie bialgebroid structure
on the base manifold.

\vspace{6pt}
{\bf Acknowledgement.}
We thank P. Etingof for discussions, which stimulated
this work.
One of the authors (A. M.) is grateful to the Bar Ilan University for hospitality, the friendly atmosphere,
and excellent research conditions.
\section{Dynamical categories}
\label{secDC}
\subsection{Hopf algebras and the double}
\label{subsecHAD}
In this subsection we fix some notation and set up general conventions concerning Hopf
algebras\footnote{For a guide in
the Hopf algebras and quantum groups the reader is referred to Drinfeld's
report \cite{Dr1} or to one of the  textbooks, e.g. \cite{K} or \cite{Mj}}
that will be used in the paper.

Let $k$ denote a field of zero characteristic
or a topological algebra of formal power series in one variable with coefficients in the field.
By an algebra we mean an associative unital algebra over  $k$; all the homomorphisms
of algebras are unital. Unless otherwise specified,
ideals are assumed to be two-sided ideals.
The symbol $\tp$ stands for the (completed) tensor product in the category of (complete) $k$-modules.

Let $\Ha$ be a Hopf algebra
over $k$ with invertible antipode $\gm$. We use the symbolic Sweedler notation
for the coproduct $\Delta(h)=h^{(1)}\tp h^{(2)}\in \Ha\tp \Ha$ and mark the tensor
components in the standard way, e.g., $\Ru=\Ru_1\tp \Ru_2\in \Ha\tp \Ha$.
We use analogous notation for an $\Ha$-coaction $\delta$ on a (left) comodule $A$,
namely,
$\delta(a)=a^{(1)}\tp a^{[2]}$, where the square brackets label
the $A$-component and the parentheses mark the component belonging to
$\Ha$. The Hopf algebra with the opposite multiplication will be
denoted by $\Ha_{op}$ whereas the Hopf algebra  with the opposite comultiplication will be denoted by
$\Ha^{op}$.

All $\Ha$-modules are assumed to be left.
Recall that an associative algebra and $\Ha$-module $\A$ is called an $\Ha$-module algebra,
or simply $\Ha$-algebra, if the action is non-degenerate (the unit acts as the identity operator),
the multiplication in $\A$ is $\Ha$-equivariant, and the unit in $\A$ generates
the trivial submodule. Recall also that a (left) $\Ha$-comodule algebra $\A$ is
an algebra and $\Ha$-comodule such that the coaction $\A\to \Ha\tp \A$ is
an algebra homomorphism.

Assuming $\Ha$ is quasitriangular with the universal R-matrix $\Ru$,
we will use the standard notation
\be
\label{eq+-}
\Ru^+=\Ru, \quad\Ru^-=\Ru^{-1}_{21}
.
\ee
The matrix $\Ru^-$ is
an alternative quasitriangular structure on $\Ha$.
We will use the following well known equalities relating the R-matrix and the antipode
\be
(\gm\tp \id)(\Ru)=\Ru^{-1}=(\id\tp\gm^{-1})(\Ru), \quad (\gm\tp \gm)(\Ru)=\Ru.
\label{R-gm}
\ee
If an  $\Ha$-algebra $\A$ satisfies the condition
\be
\label{eqQCDYA}
\la \mu= \bigl(\Ru_2 \tr \mu\bigr) \:\bigl(\Ru_1\tr\la\bigr),
\ee for all $\la,\mu\in \A$,
then $\A$ is called $\Ru$-commutative or $\Ha$-commutative (or simply quasi-commutative
if $\Ha$ and $\Ru$ are clear from the context). Note that this definition is independent on the choice of the
matrices $\Ru^\pm$.

By the dual $\Ha^*$ to Hopf algebra $\Ha$ we understand a Hopf algebra equipped with the non-degenerate Hopf
pairing
$\langle ., .\rangle \colon\Ha\tp \Ha^*\to k$.

Twist by a cocycle  $\F\in \Ha\tp \Ha$ of a Hopf algebra $\Ha$ with the
coproduct $\Delta$ is a Hopf algebra with the same multiplication
and with the coproduct $h\mapsto \F^{-1}\Delta(h)\F$.
Given two Hopf algebras $\A$ and $\B$, a bicharacter $\F$ is a non-zero element from $\B\tp \A$
obeying
$$
(\Delta_\B\tp \id)(\F)=\F_{13}\F_{23}\in \B\tp \B\tp \A,
\quad
(\id\tp \Delta_\A)(\F)=\F_{13}\F_{12}\in \B\tp \A\tp \A.
$$
A bicharacter defines a cocycle in the Hopf algebra $\A\tp \B$, see \cite{RS};
the corresponding twisted Hopf algebra \tw{\A}{\F}{\B} is
called a twisted tensor product of $\A$ and $\B$. The comultiplication in
\tw{\A}{\F}{\B} has the form
\be\label{twisttp}
\Delta(a\tp b)=(a^{(1)}\tp \F^{-1}_1 b^{(1)}\F_1)\tp (\F^{-1}_2 a^{(2)}\F_2 \tp b^{(2)})
\ee
It is convenient for our exposition to define the double $\D\Ha$ of the Hopf algebra
$\Ha$, \cite{Dr1}, as
a double cross product $\Ha\bowtie \Ha^*_{op}$, \cite{Mj}.
This is equivalent to the standard definition of the double
as $\Ha\bowtie \Ha^{*op}$, having in mind the isomorphism between $\Ha^*_{op}$
and $\Ha^{* op}$ realized via the antipode.
Algebraically,  $\D\Ha$ is dual to the tensor product
$\Ha^*\tp \Ha^{op}$ twisted by the canonical element $\sum_i e_i\tp e^i\in \Ha^{op}\tp \Ha^*$
of the pairing $\langle.,.\rangle$, where $\{e_i\}$ is the basis
in $\Ha$ and  $\{e^i\}$ its dual in  $\Ha^*_{op}$.
Explicitly, the cross relations between elements of $\Ha$ and $\Ha^*_{op}$
are given by
$$
\eta^{(1)}* h^{(1)}\langle \eta^{(2)}, h^{(2)}\rangle=\langle \eta^{(1)}, h^{(1)}\rangle h^{(2)}*\eta^{(2)},
$$
$h\in \Ha$, $\eta\in \Ha^*_{op}$. Then $\Theta=\sum_i e^i\tp e_i$ is naturally considered
as an element from the tensor square of $\D\Ha$ is a universal R-matrix of $\D\Ha$.

We will also deal with the situation when $\Ha$ is a Hopf subalgebra in another Hopf algebra, $\U$.
Then we can define a generalized double $\Ha\bowtie \U^*_{op}$ as the dual to
the twisted tensor product of $\Ha^*$ and $\U^{op}$
(this twist is induced from the subalgebra $\Ha^*\tp \Ha^{op}\subset\Ha^*\tp \U^{op}$).
Clearly the projection $\U^*\to \Ha^*$ extends to a Hopf algebra map $\Ha\bowtie \U^*_{op}\to \D\Ha$.

A quasitriangular structure $\Ru$ on $\Ha$ defines two Hopf algebra maps $\Ru^\pm\colon \Ha^*_{op}\to \Ha$
given by
\be
\Ru^\pm(\eta)=\langle\Ru^\pm_2,\eta\rangle \Ru^\pm_1,
 \quad \eta \in \Ha^*_{op}.
\label{+-maps}
\ee
These maps extend to Hopf algebra epimorphisms $\D\Ha\to \Ha$,
\be
\label{eqmaps+-}
x\tp\eta\mapsto x\Ru^+(\eta), \quad
x\tp\eta\mapsto x\Ru^-(\eta), \quad x\tp \eta \in \Ha\bowtie \Ha^*_{op}.
\ee
The universal R-matrix $\Theta$ of the double goes over into $\Ru^\pm$ under (\ref{eqmaps+-}).

\subsection{Base algebras}
\label{susecBA}
Recall that  one can assign to any monoidal category $\Cc$  a braided category  $Z(\Cc)$ called the center of $\Cc$.
Its objects are the pairs $(X,\si)$, were $X$ is an object of $\Cc$ equipped with
a family of natural isomorphisms $\si=\{\si_A\}$, $X\tp A\stackrel{\si_{A}}{\longrightarrow}A\tp X$ for all objects of $\Cc$
(these permutations should satisfy certain functorial conditions, see \cite{K}).
When $\Cc$ is a category of $\Ha$-modules, $Z(\Cc)$ is equivalent to
the category of modules over the double $\D\Ha$.
\begin{definition}{\cite{DM1}}
\label{defBA}
Let $\Cc$ be a monoidal category and $Z(\Cc)$ its center.
A commutative algebra in $Z(\Cc)$ is called a $\Cc$-\select{base algebra}.
\end{definition}
\noindent
When $\Cc$ is a category of $\Ha$-modules, we use the
term $\Ha$-\select{base algebra}.  An $\Ha$-base
algebra can be alternatively defined as a $\D\Ha$-commutative algebra.

Equivalently, an $\Ha$-base algebra can be defined as an $\Ha$-module algebra and
simultaneously a left $\Ha$-comodule algebra satisfying the conditions
\be
\label{ba1}
\delta(h\tr\la)&=&h^{(1)}\la^{(1)}\gm(h^{(3)})\tp  h^{(2)}\tr \la^{[2]},
\\
\label{ba2}
\lambda\mu&=&(\lambda^{(1)}\rhd\mu)\;\lambda^{[2]},
\ee
for all $\la,\mu\in \Lc$ and $h\in \Ha$.
This definition is equivalent to the definition of base algebra
given in \cite{DM1}.
Remark that an $\Ha$-module
and $\Ha$-comodule fulfilling condition (\ref{ba1}) is called
a Yetter-Drinfeld module.

Let $\Theta$ denote the standard quasitriangular structure on $\D\Ha$.
For simplicity, we think of $\Theta$ as an
element of $\D\Ha\tp\D\Ha$. For $\Ha$ finite dimensional, it is a canonical element
of the Hopf pairing between $\Ha$ and $\Ha^*$, $\Theta\in \Ha^*_{op}\tp \Ha$. Such an interpretation is valid
for infinite dimensional Hopf algebras close to universal enveloping algebras
and their quantizations, if $\Ha^*$ is understood as a restricted dual, and
the tensor product is completed in some topology.
In terms of the R-matrix $\Theta=\Theta_1\tp\Theta_2\in (\D\Ha)^{\tp 2}$ of the double, the
coaction $\delta$ reads
\be
\label{coaction}
\la\mapsto \lambda^{(1)}\tp\lambda^{[2]}=\Theta_2\tp \Theta_1\tr \la.
\ee
Actually, in our constructions we may understand by $\D\Ha$    \select{any}
quasitriangular Hopf algebra that contains $\Ha$ and whose universal R-matrix belongs to $\D\Ha\tp \Ha$.
Then any $\D\Ha$-commutative algebra belongs to  the center of the category
of $\Ha$-modules and therefore is an $\Ha$-base algebra. The $\Ha$-coaction is
expressed by the formula (\ref{coaction}) with $\Theta$ replaced by $\Ru$.

\begin{remarks}
\label{rem1}
Let $\Lc$ be an $\Ha$-base algebra. Then we can state the following.
\begin{enumerate}
\item
$\Lc$ is also an $\Ha^*_{op}$-base algebra, as readily follows from the definition.
The corresponding $\Ha^*_{op}$-coaction is given by $\la\mapsto \Theta^-_2\tp \Theta^-_1\tr \la$,
see notation (\ref{eq+-}).
\item
If $\Ha$ is quasitriangular
and $\Lc$ is $\Ha$-commutative, then $\Lc$ has two $\Ha$-base algebra
structures defined by  $\Ru^\pm$, where $\Ru$ is the
R-matrix of $\Ha$. Namely, the double $\D\Ha$ acts on
$\Lc$ through  the  projections (\ref{eqmaps+-}) to $\Ha$.
The Hopf algebra homomorphisms (\ref{eqmaps+-}) sends
$\Theta^\pm$ to $\Ru^\pm$, hence the algebra $\Lc$ is $\D\Ha$-commutative.
In terms of the R-matrix, the $\Ha$-coactions are given by
\be
\label{d+-}
\delta^+(\la)= \Ru^-_2\tp \Ru^-_1\tr \la =\Ru^{-1}_1\tp \Ru^{-1}_2\tr \la
,\quad
\delta^-(\la)= \Ru^+_2\tp \Ru^+_1\tr \la=\Ru_2\tp \Ru_1\tr \la
.
\ee
We denote by $\Lc_\pm$ the two $\Ha$-base algebra
structures on $\Lc$ that correspond to the coactions
$\delta^\pm$.
\item
Combining two previous remarks, we state that an $\Ha$-base algebra has two
{\em different} $\D\Ha$-base algebra structures. The $\Ha$- and $\Ha^*_{op}$-coactions
expressed through $\Theta^\pm$ may be considered as $\D\Ha$-coactions
via the  embeddings of $\Ha$ and $\Ha^*_{op}$ into $\D\Ha$.
The $\D\Ha$-coactions are given by
$\la\mapsto \Theta^\pm_2\tp \Theta^\pm_1\tr \la$.
\item
Let us fix that $\D\Ha$-base algebra structure on $\Lc$ which
 corresponds to the $\Ha$-coaction, cf. the previous remark.
Assume that $\Ha$ is a Hopf subalgebra in a Hopf algebra $\U$, thus there is
a natural projection $\U^*\to \Ha^*$ inducing an $\U^*_{op}$-action on $\Lc$.
Then $\Lc$ is a base algebra over the generalized double $\Ha\bowtie \U^*_{op}$.
\end{enumerate}
\end{remarks}

\begin{lemma}\label{lemcenter}
 Let $\Lc$ be an $\Ha$-base algebra. Then any $\Ha$-invariant element in $\Lc$
 belongs to the center $Z(\Lc)$.
\end{lemma}

\begin{proof}
Let $\mu\in\Lc$ be $\Ha$-invariant. Then for any $\la\in\Lc$
one has $\lambda\mu=(\lambda^{(1)}\rhd\mu)\:\lambda^{[2]}=
\varepsilon(\lambda^{(1)})\:\mu\:\lambda^{[2]}=\mu\lambda$. The first equality
follows from the $\D\Ha$-commutativity of $\Lc$.
\end{proof}

\begin{definition}
\label{defnormalBA}
An $\Ha$-base algebra $\Lc$ is called quasi-transitive  if $\Lc^{\D\Ha}$, the set of $\D\Ha$-invariant
elements in $\Lc$, coincides with $k$.
\end{definition}

It follows from Lemma \ref{lemcenter} that $\Lc^{\D\Ha}$ is a commutative
algebra belonging to the center of $\Lc$. Let  $\chi$ be a character of $\Lc^{\D\Ha}$, i.e.
a one dimensional representation.
Consider the ideal $J_\chi$ in $\Lc$ generated by the kernel of  $\chi$.
\begin{propn}
\label{normal}
The quotient $\Lc/J_\chi$ is a quasi-transitive $\Ha$-base algebra.
\end{propn}
\begin{proof}
The ideal $J_\chi$ is obviously $\D\Ha$-invariant, hence the quotient $\Lc/J_\chi$
is an $\D\Ha$-algebra. It is quasi-commutative, being a quotient of a quasi-commutative
algebra. By construction, the subalgebra of invariants in $\Lc/J_\chi$ coincides with $k$.
\end{proof}
\begin{examples}

Let us give some examples of base algebras. A detailed consideration to some of them is given in \cite{DM1}.
\begin{enumerate}
\item
$\Ha$ itself is an (quasi-transitive) $\Ha$-base algebra, being equipped with the  adjoint action and the coproduct
coaction.
\item
$\Ha^*_{op}$ is a (quasi-transitive) $\Ha$-base algebra due to the symmetry
$\Ha \leftrightarrow \Ha^*_{op}$ in the definition of base algebras.
\item
Consider an FRT algebra associated with a finite dimensional representation of
a quasitriangular Hopf algebra $\Ha$, \cite{FRT}. It is a commutative algebra
in the category of $\Ha$-bimodules, whence it is an $\Ha\tp \Ha_{op}$-base algebra (cf. Remark
\ref{rem1}.2).
\item
Assume again that $\Ha$ is quasitriangular. A reflection equation algebra
(studied in \cite{KSkl}) is, in fact, a commutative algebra in the category of modules
over the twisted tensor product \tw{\Ha}{\Ru}{\Ha}, \cite{DM2}. Therefore it is
$\D\Ha$-commutative and thus an $\Ha$-base algebra.
\item
Let $\Lc$ and $\Lc_1$ be two $\Ha$-base algebras. On the linear space
$\Lc\tp\Lc_1 $ define an associative algebra structure by the multiplication
$$
(\la\tp \mu)(\al\tp \bt)\!:=\la(\Theta_2\tr \al)\tp (\Theta_1\tr \mu)\bt.
$$
This algebra is a braided tensor product
of two $\D\Ha$-commutative algebras,  hence it is \tw{\D\Ha}{\Theta}{\D\Ha}-commutative and
has two structures of $\D\Ha$-base algebra.
In case $\Lc_1=\Ha$, it coincides with the smash product $\Lc\rtimes \Ha$ as an associative algebra.

\end{enumerate}
\end{examples}

\subsection{Dynamical categories}
The notion of dynamical extension (dynamization) of a monoidal category admits various formulations, \cite{DM1}, which
become equivalent under  certain circumstances. We will work with a
category $\Mc_\Ha$ of $\Ha$-modules and its extension  $\overline{\Mc}_{\Ha,\Lc}$
over an $\Ha$-base algebra $\Lc$ in the sense of the following definition.
\begin{definition}{\cite{DM1}}
\label{defDC}
\select{Dynamization} of the category $\Mc_\Ha$ over the $\Ha$-base algebra
$\Lc$ is a strict monoidal category $\overline{\Mc}_{\Ha,\Lc}$ defined by
the following conditions
\begin{enumerate}
\item
objects of $\overline{\Mc}_{\Ha,\Lc}$ are the objects of $\Mc_\Ha$,
\item
$\Hom_{\overline{\Mc}_{\Ha,\Lc}}(X,Y)$ is the set of $\Ha$-equivariant linear maps
from $X$ to $Y\tp \Lc$. The composition $\phi\circ \psi$ of
two morphisms $\phi\in \Hom_{\overline{\Mc}_{\Ha,\Lc}}(X,Y)$ and
$\psi\in \Hom_{\overline{\Mc}_{\Ha,\Lc}}(Y,Z)$ is the composition map
$$
X\stackrel{\phi}{\longrightarrow} Y\tp \Lc\stackrel{\psi\tp \id_\Lc}{\longrightarrow} Z\tp \Lc\tp \Lc
\stackrel{\id_Z\tp \mathrm{m}_\Lc}{\longrightarrow} Z\tp \Lc,
$$
where $\mathrm{m}_\Lc$ is the multiplication in $\Lc$,
\item
tensor product of objects from $\overline{\Mc}_{\Ha,\Lc}$ is the same as in $\Mc_\Ha$,
\item
tensor product of  morphisms $\phi\in \Hom_{\overline{\Mc}_{\Ha,\Lc}}(X,X')$ and
$\psi\in \Hom_{\overline{\Mc}_{\Ha,\Lc}}(Y,Y')$
is given by the composition
$$X\tp Y\stackrel{\phi\tp \psi}{\longrightarrow} X'\tp \Lc\tp Y'\tp \Lc
\stackrel{\tau_{Y'}}{\longrightarrow}X'\tp Y'\tp\Lc\tp \Lc
\stackrel{\mathrm{m}_\Lc}{\longrightarrow}X'\tp Y'\tp\Lc.
$$
where $\tau_{Y'}$ is the permutation $\Lc\tp Y'\to Y'\tp \Lc$ expressed via the $\Ha$-coaction
on $\Lc$
by the formula $\la\tp y\mapsto \la^{(1)}\tr y\tp \la^{[2]}$.
\end{enumerate}
\end{definition}
The category $\overline{\Mc}_{\Ha,\Lc}$ generalizes the category of Etingof
and Varchenko that was constructed
in  \cite{EV2} for commutative cocommutative $\Ha$ and $\Lc$ being a certain extension of $\Ha$.
The category $\overline{\Mc}_{\Ha,\Lc}$ was introduced in \cite{DM1} in order to formulate
the classical and quantum dynamical Yang-Baxter equations for an arbitrary Lie bialgebras and their quantizations.
The purpose of the present paper is to realize $\overline{\Mc}_{\Ha,\Lc}$ and its important subcategories
via representations of bialgebroids.
The notion of bialgebroid is a generalization of the notion of Hopf algebra, \cite{Lu}.
The next section is a brief introduction to this theory.

\section{Some basics on bialgebroids}
\label{secHA}
\subsection{General definition and examples}
The reconstruction theorem states that a fiber functor from a monoidal category $\Cc$
to the monoidal category of vector spaces gives rise to a bialgebra whose category of representations
is equivalent to $\Cc$, see e.g. \cite{Mj}. Not all monoidal categories admit such a fiber functor,
thus not all of
them are related to bialgebras, \cite{GK}. A more general concept of
 functor to the monoidal category of bimodules over some associative algebra
leads to the notion of bialgebroid \cite{Lu}. Similarly to the bialgebra case,
representations of a bialgebroid also form a monoidal category.
\begin{definition}
\label{Lu}
Let $\Rg$ be an associative unital algebra over $k$. An associative
unital algebra $\Bg$ over $k$ is called
a \select{ bialgebroid over  base} $\Rg$ or $\Rg$-bialgebroid if there exist
\begin{enumerate}
\item an algebra homomorphism $s\colon \Rg\to \Bg$
(source map) and an algebra anti-homomorphism $ t \colon \Rg\to \Bg$
(target map) making $\Bg$ an $\Rg$-bimodule by
$\la\llcorner a \!:=  s (\la) a$, $a\lrcorner\la \!:=
 t (\la)a$, $\la\in \Rg$, $a\in \Bg$,
\item a coassociative bimodule map
(comultiplication) $\Delta\colon \Bg\to \Bg\tp _{\Rg} \Bg$ which is a
homomorphism into the unital associative algebra specified by the
condition
\be
\{z\in \Bg\tp_\Rg \Bg| \>z\bigl( t (\la)\tp 1\bigr) = z\bigl(1\tp  s (\la)\bigr), \forall \la\in \Rg\},
\label{subalg}
\ee
\item a bimodule map
(counit) $\ve\colon \Bg\to \Rg$ such that $\ve(1_\Bg)=1_\Rg$,
\be
&\ve\bigl(a\; (s\circ\ve) (b)\bigr)=\ve(ab)=\ve\bigl(a\; (t\circ\ve) (b)\bigr), \quad a,b\in \Bg,
\quad \mbox{ and}
\label{anchor}
\\[7pt]
&(\ve\tp_\Rg\id_\Bg)\circ \Delta =\id_\Bg= (\id_\Bg\tp_\Rg\ve)\circ
\Delta \label{counit}
\ee
under the identification
$\Rg\tp_\Rg\Bg\simeq \Bg\simeq  \Bg\tp_\Rg\Rg$.
\end{enumerate}
\end{definition}

\begin{remarks}

\begin{enumerate}
\item
The images of the source and target maps  in $\Bg$ commute, by virtue of condition 1.
\item
In general, the tensor product $\Bg\tp_\Lc\Bg$ has no natural structure of associative algebra. However, the
element $z(a\tp b):=z_1a\tp_\Lc z_2b\in \Bg\tp_\Lc\Bg$ is well defined for any $z\in \Bg\tp_\Lc \Bg$
and $a\tp b\in \Bg\tp\Bg$. The condition 2 selects a natural algebra in $\Bg\tp_\Lc\Bg$.
\item
Since $\Delta$ is a  bimodule map, one has $\Delta\circ s= s\tp_\Rg 1$ and $\Delta\circ t= 1\tp_\Rg t$.
\item
Condition 3 implies the identities $\ve\circ s=\ve\circ t=\id_\Rg$ and
makes $\Rg$ a left $\Bg$-module by
\be
a\vdash\la:=\ve\bigl(a s(\la)\bigr)=\ve\bigl(a t(\la)\bigr),
\quad
a\in \Bg, \la\in \Rg,
\label{regular}
\ee
where the right equality is a consequence of (\ref{anchor}).
The $\Rg$-bimodule structure on $\Rg$ induced by this action coincides with
the standard one. The action (\ref{regular}) is called \select{anchor}.
One can check that
\be
a\: s(\la) = s(a^{(1)}\vdash\la)\:a^{(2)},
\quad
a\: t(\la) = t(a^{(2)}\vdash\la)\:a^{(1)}.
\label{s-anch-t}
\ee
Sometimes the anchor is  introduced
separately; then the condition (\ref{anchor}) is dropped from
definition of bialgebroid, see \cite{Lu}. In our definition we follow \cite{Szl}.
\end{enumerate}
\end{remarks}

Any left $\Bg$-module $V$ is a natural $\Rg$-bimodule. We call this correspondence
the forgetful functor.
Given two $\Bg$-modules $V$ and $W$,
the tensor product $V\tp_\Rg W$ acquires a left $\Bg$-module structure
via the coproduct, due to condition 2 of Definition \ref{Lu}.
The whole set of axioms from Definition \ref{Lu} ensures that the
left $\Bg$-modules form a monoidal category,
with $\Rg$ being the unit object. The forgetful functor  to
the category of $\Rg$-bimodules is strong monoidal, i.e. preserves tensor products.
Conversely,  suppose a pair of algebras $(\Bg,\Rg)$ satisfies condition 1
of Definition \ref{Lu}
and there is a  monoidal structure on the category of left $\Bg$-modules.
Suppose the forgetful functor
to the category of $\Rg$-bimodules is strong monoidal.
Then $\Bg$ is an $\Rg$-bialgebroid, see \cite{Szl}.
\begin{remark}
\label{right}
The bialgebroid $\Bg$ from Definition  \ref{Lu} is a \select{left} one. This means that
the $\Rg$-bimodule structure on $\Bg$ is defined by the source and target maps and multiplication from the left.
Alternatively, one can consider $\Bg$ as a $\Rg$-bimodule using multiplication from the right and
require that the right $\Bg$-modules form a monoidal category with the forgetful functor to a $\Rg$-bimodules.
Such bialgebroids are called \select{right} ones; one can readily recover their definition by the apparent modification of
Definition  \ref{Lu}. Although right modules over $\Bg$ are the same as left modules over $\Bg_{op}$,
sometimes the notion of right bialgebroid proves to be convenient to work with.
\end{remark}
\begin{example}[Bialgebroid $\End(\Rg)$]
Let $\Rg$ be a finite dimensional associative unital algebra over the field $k$. Denote by
$\E$ the algebra of endomorphisms of $\Rg$ over $k$.
For $a\in \Rg$ let $L_a$ and $R_a$ be the linear operators acting on $\Rg$ via the  left and right multiplication
by $a$; they define an algebra and anti-algebra maps from $\Rg$ to $\E$, respectively.
Thus $\E$ is a natural $\Rg$-bimodule: the element $a\tp b\in \Rg\tp_k \Rg_{op}$ acts on $\E$ by
multiplication by $L_a R_b$ from the left.
The algebra $\E$ is in fact an $\Lc$-bialgebroid with the coproduct defined by
$\Delta(f)(a\tp b):= f(ab)$ and the counit $\ve(f):=f(e)$, see \cite{Lu}.
\end{example}
\begin{example}[Bialgebroid structure on $\Rg\tp \Rg_{op}\tp \Ha$]
\label{R-R-H}
Suppose $\Rg$ is a left $\Ha$-module algebra for some Hopf algebra $\Ha$.
The action of $\Ha$ on $\Rg$ is denoted by $\tr$.
Consider the associative algebra $\Bg$ built on $\Rg\tp \Rg_{op}\tp \Ha$
and equipped with the multiplication
$$
(\la\tp \mu \tp f)(\zeta\tp \eta \tp g):=\la (f^{(1)}\tr \zeta)\tp \mu (f^{(3)}\tr \eta)\tp f^{(2)}g.
$$
Let $\iota$ denote the (anti-algebra) identity map from $\Rg$ to $\Rg_{op}$.
It is not difficult  to show that $\Bg$ is a bialgebroid with
the source map $ s  \colon \la\mapsto \la\tp 1\tp 1$
the target map $ t \colon \la\mapsto 1\tp \iota(\la)\tp 1$, the coproduct
$\Delta(\la\tp \mu\tp h):=\bigl(\la\tp 1\tp h^{(1)}\bigr)\tp_\Lc \bigl(1\tp \mu\tp h^{(2)}\bigr)$
and the counit $\ve(\la\tp \mu\tp h):=\la\;\iota^{-1}(\mu)\ve(h)$.
The anchor action (\ref{regular}) is given explicitly by
$$
(\la\tp\mu\tp h)\llcorner\zeta= \la(h\tr \zeta)\iota^{-1}(\mu),
$$
for $(\la\tp\mu\tp h)\in\Bg$ and $\zeta \in \Rg$.
\end{example}
\begin{example}[Bialgebras]
A bialgebra over the field $k$ is a bialgebroid whose base is $k$.
\end{example}
\begin{example}[Tensor product of bialgebroids]
\label{H-B}
Let $(\Bg_i,\Rg_i,s_i,t_i,\Delta_i,\ve_i)$, $i=1,2$,
be a pair of bialgebroids. Then one
can build their tensor product bialgebroid over the base $\Rg_1\tp \Rg_2$.
As an associative algebra, this is the standard tensor product $\Bg_1\tp \Bg_2$.
The source, target, and counit maps are respectively $s_1\tp s_2$, $t_1\tp t_2$,
and $\ve_1\tp \ve_2$.
The coproduct is given
by
$$\Delta(x\tp y):= (x^{(1)}\tp y^{(1)})\tp_{(\Rg_1\tp\Rg_2)}  (x^{(2)}\tp y^{(2)}).$$
In particular, if one of the bialgebroids, say $\Bg_1$ is a Hopf algebra,
then the tensor product bialgebroid  will be over the base $\Rg_2$.
\end{example}
\begin{definition}
Let $(\Bg_i,\Lc,s_i,t_i,\Delta_i,\ve_i)$, $i=1,2$, be two $\Rg$-bialgebroids.
An algebra map
$\varphi\colon \Bg_1\to \Bg_2$ is called a homomorphism of  bialgebroids
if it is an $\Rg$-bimodule map and
\be
\label{bial_hom}
\ve_2\circ \varphi =\ve_1,\quad (\varphi\tp_\Rg \varphi)\circ \Delta_1 =\Delta_2\circ \varphi.
\ee
\end{definition}
\begin{example}
For an arbitrary bialgebroid $\Bg$ over a finite dimensional base $\Rg$ the anchor
map  $\Rg\to \End_k(\Rg)$ is a homomorphism of Lie bialgebroids.
\end{example}
Of particular interest for us will be the notion of quotient bialgebroid.
\begin{definition}
Let $\Bg$ be an $\Rg$-bialgebroid. A two-sided ideal $J$ in the algebra $\Bg$
is called a \select{biideal} if $\Delta(J)\subset J\tp_\Rg \Bg +\Bg\tp_\Rg J$
and $\ve(J)= 0$.
\end{definition}
\noindent
Given a biideal $J\subset \Bg$ the quotient $\Bg/J$ is naturally endowed with
an $\Rg$-bialgebroid structure such that the projection $j\colon\Bg\to \Bg/J$
is a bialgebroid homomorphism.
\begin{remark}
Note that any biideal lies in the kernel of  the anchor map
since the latter is expressed through the counit by formula (\ref{regular}).
\end{remark}
\subsection{Quasitriangular structure and twist}
\label{subsecQST}
In the Hopf algebra theory, a quasitriangular structure
on a Hopf algebra is essentially the same as a braiding in
the monoidal category of its modules.
Analogously to Hopf algebras one can define quasitriangular bialgebroids, with inevitable
complications caused by non-commutativity of the base. A quasitriangular
structure on a bialgebroid gives rise to a braiding in the category of
its modules.

Let $\Bg$ be an $\Rg$-bialgebroid.
Then every $\Bg$-module,
and $\Bg$ in particular, is
also a natural $\Rg_{op}$-bimodule
with respect to the left and right $\Rg_{op}$-actions defined through the
target and source maps, correspondingly.
Given two $\Bg$-modules $M_1$ and $M_2$, the flip
$M_1\tp M_2\to M_2\tp M_1$
induces an invertible map $\sigma_{M_1,M_2}\colon M_1\tp_\Rg M_2\to M_2\tp_{\Rg_{op}}M_1$.
Let us define a structure of an $\Rg_{op}$-bialgebroid, $\Bg^{op}$, on
the algebra $\Bg$.
The target and source maps from $\Rg$ to $\Bg$ viewed
as algebra and anti-algebra maps from $\Rg_{op}$ to $\Bg$
give, respectively, the source and target maps
of the $\Rg_{op}$-bialgebroid $\Bg^{op}$.
To define $\Bg^{op}$, it is enough to specify the corresponding
monoidal structure on the left $\Bg$-modules.
Let us define a new tensor product of two $\Bg$-modules $M_1$ and $M_2$
as the $\Rg_{op}$-bimodule $M_1\tp_{\Rg_{op}} M_2$ equipped  with  the following $\Bg$-action:
\be
\label{coop}
M_1\tp_{\Rg_{op}} M_2\stackrel{\si^{-1}_{M_2,M_1}}{\longrightarrow} M_2\tp_{\Rg} M_1
\stackrel{\Delta(a)}{\longrightarrow} M_2\tp_{\Rg} M_1
\stackrel{\si_{M_2,M_1}}{\longrightarrow}M_1\tp_{\Rg_{op}} M_2
,\quad a\in \Bg
\ee
This tensor product is associative, as follows
from the coassociativity of $\Delta$ and the "hexagon" identity obeyed by $\si$.
One can check that the corresponding coproduct $\Bg^{op}$ is given by
$\Delta^{op}=\si_{\Bg,\Bg}\circ\Delta$, and the counit is
$\iota\circ\ve_\Bg$, where $\iota$ is the anti-isomorphism $\Rg\to \Rg_{op}$ implemented
by the identity map.
\begin{definition}
\label{R-matrix}
A bialgebroid $\Bg$ is called \select{quasitriangular} if
there is a monoidal isomorphism  $\Mod \Bg\to \Mod \Bg^{op}$
identical on objects, and the transformation
of tensor products is defined by an element
$\Ru=\Ru_1\tp_{\Rg_{op}} \Ru_2\in \Bg\tp_{\Rg_{op}} \Bg$ (universal R-matrix):
$$
M_1\tp_\Rg M_2 \stackrel{\Ru}{\longrightarrow} M_1\tp_{\Rg_{op}} M_2,
\quad
x_1\tp_\Rg x_2 \mapsto \Ru_1x_1\tp_{\Rg_{op}} \Ru_2 x_2,
$$
for any pair of modules $M_1,M_2$.
\end{definition}

It follows that there exists an element $\bar\Ru\in \Bg\tp_\Rg \Bg$ such that
$
\Ru \bar\Ru = 1\tp_{\Rg_{op}}1$ and $\bar\Ru \Ru  = 1\tp_{\Rg}1.
$
The element $\bar \Ru$ implements the inverse isomorphism $\Mod \Bg^{op}\to \Mod \Bg$ and
it is a quasi-triangular  structure on the coopposite bialgebroid $\Bg^{op}$.

We will use the term \select{quantum groupoid} for a quasitriangular bialgebroid.
\begin{propn}
\label{R-matrix_prop}
An element $\Ru\in \Bg\tp_{\Rg_{op}} \Bg$ defines a quasitriangular structure
on $\Bg$ if and only if
\begin{enumerate}
\item
$\Ru\bigl( t (\la)\tp 1\bigr)=\Ru\bigl(1\tp s (\la)\bigr)$, for all $\la\in \Rg$,
\item
 for all $a\in \Bg$ equation $\Ru\Delta(a)=\Delta^{op}(a)\Ru$ holds in $\Bg\tp_{\Rg_{op}} \Bg$,
\item
 equations
\be
(\Delta^{op}\tp_{\Rg_{op}}\id)(\Ru)&=&\Ru_{23}\Ru_{13}:= \Ru_1\tp_{\Rg_{op}} \Ru(1\tp \Ru_2),
\label{hex1}
\\
(\id\tp_{\Rg_{op}}\Delta^{op})(\Ru)&=&\Ru_{12}\Ru_{13}:=\Ru(\Ru_1\tp 1)\tp_{\Rg_{op}}  \Ru_2)
\label{hex2}
\ee
hold in $\Bg\tp_{\Rg_{op}} \Bg\tp_{\Rg_{op}} \Bg$,
\item
there exists an element $\bar\Ru\in \Bg\tp_\Rg \Bg$ such that
$\bar\Ru\bigl( s (\la)\tp 1\bigr)=\bar\Ru\bigl(1\tp t (\la)\bigr)$, for all $\la\in \Rg$,
and
$
\Ru \bar\Ru = 1\tp_{\Rg_{op}}1$ and $\bar\Ru \Ru  = 1\tp_{\Rg}1.
$
\end{enumerate}
\end{propn}
\begin{proof}
A direct computation.
\end{proof}
\begin{remarks}
Let us make a few comments on the conditions of Proposition \ref{R-matrix_prop}.
\begin{enumerate}
\item
By condition 1,
one has $\Ru(a\tp_\Rg b )\in \Bg\tp_{\Rg_{op}} \Bg $
for all $a,b\in \Bg$. Analogously, condition 4 implies
$\bar\Ru(a\tp_{\Rg_{op}} b )\in \Bg\tp_\Rg \Bg $
for all $a,b\in \Bg$.
The R-matrix $\Ru$ lives, in fact, in $\Bg^{op}\tp_{\Rg_{op}} \Bg^{op}$.
This explains appearance
of the opposite coproduct in (\ref{hex1}-\ref{hex2}). On the contrary,
the inverse $\bar\Ru$ is supported in $\Bg\tp_{\Rg} \Bg$.
\item
Both sides of the equation from condition 2 are well defined, cf. remark 1.
\item
The right-hand side expressions in (\ref{hex1}-\ref{hex2}) are
correctly defined, i.e. are independent on the representative
$\Ru_1\tp \Ru_2$ of $\Ru\in \Bg\tp_{\Rg_{op}} \Bg$  in $\Bg\tp \Bg$.
Indeed, from  condition 2 one deduces
\be
\Ru\bigl(1\tp t (\la)\bigr)=t (\la)\Ru_1\tp_{\Lc_{op}} \Ru_2
,\quad
\Ru\bigl( s (\la)\tp1\bigr)=\Ru_1\tp_{\Lc_{op}} s (\la)\Ru_2
\label{permutation}
\ee
for all $\la\in \Rg$.
Equations (\ref{permutation}) imply that the two maps $\Bg\tp \Bg\to \Bg\tp_{\Rg_{op}}\Bg\tp_{\Rg_{op}}\Bg$
defined by
\be
\label{auxi_maps}
\hat\jmath_{23}\colon x\tp y\mapsto x \tp_{\Rg_{op}}\Ru(1\tp y),
\quad
\hat\jmath_{12}\colon x\tp y\mapsto \Ru(x\tp 1)  \tp_{\Rg_{op}}y,
\ee
are factored through maps
$\jmath_{12},\jmath_{23}\colon\Bg\tp_{\Rg_{op}} \Bg\to \Bg\tp_{\Rg_{op}}\Bg\tp_{\Rg_{op}}\Bg$.
Let us check this, say, for $\hat\jmath_{23}$.
In view of (\ref{permutation}), we have
$$
x \tp_{\Rg_{op}} \Ru\bigl(1\tp  t (\la)y\bigr)=
x \tp_{\Rg_{op}} t (\la)\Ru_1\tp_{\Lc_{op}} \Ru_2 y=
 s (\la)x \tp_{\Rg_{op}} \Ru(1\tp y)
$$
for all $x,y\in \Bg$, $\la\in \Rg$.
This shows that $\hat\jmath_{23}\bigl(x \tp  t (\la)y\bigr)= \hat\jmath_{23}\bigl(s (\la)x \tp y\bigr)$, i.e.
the value of $\hat\jmath_{23}$ depends only on the class of a representative of
$\Bg\tp_{\Rg_{op}} \Bg$ in $\Bg\tp\Bg$.

Now notice that the right-hand sides of equations (\ref{hex1}-\ref{hex2})
are equal to $\jmath_{23}(\Ru)$ and  $\jmath_{12}(\Ru)$, respectively.
\end{enumerate}
\end{remarks}
\begin{thm}
Suppose a bialgebroid $\Bg$ is quasitriangular.
Then the collection of morphisms
$\si^{-1}_{M_2,M_1}\circ(\rho_1\tp_{\Rg_{op}} \rho_2)(\Ru)
\in \Hom_{\Bg}(M_1\tp_{\Rg} M_2,M_2\tp_{\Rg} M_1)$, where $(M_i,\rho_i)$, $i=1,2$, are
 $\Bg$-modules, is a braiding
in the monoidal category $\Mod \Bg$.
\end{thm}
\begin{proof}
Follows from the definition of $\Ru$.
\end{proof}

Analogously to Hopf algebras, one can consider twists of bialgebroids.
\begin{definition}[\cite{Xu1}]
An element $\Psi=\Psi_1\tp_\Rg\Psi_2\in \Bg\tp_\Rg\Bg$, where
$\Bg$ is an $\Rg$-bialgebroid, is called a twisting cocycle if
\be
\Delta(\Psi_{1})\Psi\tp_\Rg \Psi_{2}=\Psi_{1}\tp_\Rg \Delta(\Psi_{2})\Psi.
\label{bialg_cocycle}
\ee
and $(\ve\tp_\Rg \id)(\Psi)=(\id \tp_\Rg \ve)(\Psi)=1 \tp_\Rg 1.$
\end{definition}
Given a twisting cocycle, the space $\Rg$ is equipped with a new
multiplication
$$
\la*\mu:=(\Psi_1\vdash \la)(\Psi_2\vdash \la),
$$
making it an associative algebra, $\tilde \Rg$.
Applying equation (\ref{bialg_cocycle}) to
$\Bg\tp\Rg\tp \Rg$, $\Rg\tp\Bg\tp \Rg$, and
$\Rg\tp\Rg\tp \Bg$, one obtains
that
\be
\tilde s(\la):= s(\Psi_1\vdash\la )\Psi_2
,\quad
\tilde t(\la):= t(\Psi_2\vdash\la )\Psi_1
,\quad \la\in \Bg,
\ee
are, respectively, an algebra and anti-algebra maps
from $\tilde \Rg$ to $\Bg$ and their images commute in $\Bg$.
Thus $\Bg$ becomes an $\tilde \Rg$-bimodule by means of the new
source and target maps, $\tilde s$ and $\tilde t$.
Applying formulas (\ref{s-anch-t}), one can check that
$\Psi\bigl(\tilde t(\la)\tp 1\bigr)=\Psi\bigl(1\tp\tilde s(\la)\bigr)$.

Thus  twisting cocycle defines an operator acting from $\Bg\tp_{\tilde\Rg}\Bg$
to $\Bg\tp_{\Rg}\Bg$ by the mapping $a\tp_{\tilde\Rg} b\mapsto \Psi_1 a\tp_{\Rg} \Psi_2 b$.
It is called invertible if there is an element $\Psi^{-1}\in \Bg\tp_{\tilde\Rg}\Bg$
such that $\Psi\Psi^{-1}\in 1\tp_{\Rg} 1$ and $\Psi^{-1}\Psi\in 1\tp_{\tilde\Rg}1$.

\begin{propn}[\cite{Xu1}]
\label{twist}
Let $\Bg$ be an $\Rg$-bialgebroid and  $\Psi\in \Bg\tp_\Rg\Bg$ be an invertible twisting cocycle.
Let $\tilde \Delta$ denote the map
\be
\label{twcopr}
a\mapsto \Psi^{-1}\Delta(a)\Psi
\ee
from
$\Bg$ to $\Bg\tp_{\tilde \Rg}\Bg$.
Then $(\Bg,\tilde \Rg,\tilde s, \tilde t, \tilde \Delta, \ve)$
is an $\tilde\Rg$-bialgebroid called the twist of $\Bg$ by $\Psi$.
\end{propn}
\begin{remark}
\label{hom_tw}
Given two $\Rg$-bialgebroids $\Bg_i$, $i=1,2$, a twist $\Psi\in \Bg_1\tp_\Rg\Bg_1$,
and a homomorphism $\varphi\colon \Bg_1\to\Bg_2$,
the element $(\varphi\tp_\Rg\varphi)(\Psi)\in \Bg_2\tp_\Rg\Bg_2$
is a twist in $\Bg_2$.
Then $\varphi$ becomes  a homomorphism of twisted bialgebroids,
$\varphi\colon \tilde \B_1\to \tilde \B_2$.
\end{remark}
A bialgebroid  twist induces a transformation of monoidal categories.
For any pair $M_1$, $M_2$ of $\Bg$-modules,
the twist $\Psi$ gives a map $M_1\tp_{\tilde \Rg}M_2\to M_1\tp_\Rg M_2$
intertwining the actions of $\tilde\Bg$ and $\Bg$.
If  $\Bg$ is quasitriangular, then the braiding in $\Mod \Bg$
defines a braiding in $\Mod \tilde\Bg$. This follows from
the following fact.
\begin{propn}
\label{Rtwisted}
Let $\Bg$ be a quasitriangular  $\Rg$-bialgebroid with the universal R-matrix $\Ru$
and let $\Psi\in \Bg\tp_{\Rg} \Bg$ be a twisting cocycle.
Then the twisted bialgebroid $\tilde \Bg$ is quasitriangular, with the universal R-matrix
$\tilde \Ru:=(\Psi_{21})^{-1}\Ru\Psi$, where $\Psi_{21}=\si_{\Bg,\Bg}(\Psi)\in \Bg\tp_{\Rg_{op}}\Bg$.
\end{propn}
\begin{proof}
First of all notice that $\tilde \Ru=(\Psi_{21})^{-1}\Ru\Psi$ is a well defined element
of $\tilde\Bg\tp_{\tilde \Rg_{op}}\tilde\Bg$.
The proof is carried out by a direct computation.
\end{proof}
\begin{remark}
\label{R-tw}
The R-matrix $\Ru$ is a special twist of the coopposite bialgebroid $\Bg^{op}$, analogously to the Hopf algebra case.
\end{remark}
\section{Bialgebroids over a quasi-commutative base}
\label{secBH_R}
\subsection{Bialgebroid $\Rg\rtimes\Ha$}
\label{subsecHR}
In this subsection we assume that the Hopf algebra $\Ha$ is quasitriangular,
with the universal R-matrix $\Ru$.
The module algebra $\Rg$ is assumed to be $\Ha$-commutative, i.e.
\be
(\Ru_2\tr \mu)(\Ru_1\tr \la)= \la \mu
\ee
for all $\la,\mu\in \Ha$.
We use the standard notation
$\Ru^+=\Ru$ and $\Ru^-=\Ru^{-1}_{21}$.
Recall that $\Ru^-$ gives an alternative quasitriangular structure on $\Ha$,
and for $\Lc$ to be quasi-commutative does not
depend on the choice of $\Ru=\Ru^\pm$.

Recall that $\Lc$ is equipped with two structures $\Lc_\pm$
of $\Ha$-base algebra corresponding to the two coactions $\delta^\pm$, cf.
Remark \ref{rem1}.
Consider the associative algebra
$\Rg\rtimes \Ha $ endowed with the smash product
multiplication
\be
(\la\tp f)(\mu \tp g)&:=&\la (f^{(1)}\tr \mu) \tp f^{(2)}g.
\ee
Introduce linear maps $ s $ and $ t^\pm$ from $\Rg$ to $\Rg\rtimes \Ha$ by
\be
\label{bt_pm}
 s(\la)= \la\otimes 1,\quad
 t ^\pm (\la)= \Ru^\pm_2\tr \la\otimes \Ru^\pm_1,
\ee
for $\la\in \Rg$.
By the construction of smash product, $ s $ is an algebra embedding.
The maps $t^\pm$ are expressed through the $\Ha$-coactions
by the formula
\be
\label{t-d}
t^\pm(\la)=\la^{[2]}\otimes \gm^{-1} (\la^{(1)}),
\ee
where $\delta^{\pm}=\la^{(1)}\tp\la^{[2]}$ and $\gm$ is the antipode in $\Ha$, cf. formulas
(\ref{R-gm}).

\begin{lemma}
\label{lemma-s-t}
The maps $ t ^\pm$ are algebra anti-homomorphisms. For every pair $\la,\mu\in
\Rg$ one has $ s (\la) \:t ^\pm(\mu)= t ^\pm(\mu) s (\la)$.
\end{lemma}
\begin{proof}
Let us check the statement for $ t = t ^+$.
Using the bicharacter properties
$$(\Delta\tp \id)(\Ru)=\Ru_{13}\Ru_{23},\quad (\id\tp \Delta)(\Ru)=\Ru_{13}\Ru_{12}$$
 we find
\be
 t (\la \mu)&=&\Ru_2\tr (\la\mu)\otimes \Ru_1=(\Ru^{(1)}_2\tr \la)(\Ru^{(2)}_2\tr \mu)\otimes \Ru_1 =
(\Ru_2\tr \la)(\Ru_{2'}\tr \mu)\otimes \Ru_{1'}\Ru_1.
\nn
\ee
On the other hand,
\be
 t (\mu) t (\la)&=&(\Ru_{2'}\tr \mu\otimes\Ru_{1'})(\Ru_2\tr \la\otimes\Ru_{1})= (\Ru_{2'}\tr \mu)
\bigl(\Ru^{(1)}_{1'}\Ru_2\tr \la\bigr)\otimes\Ru^{(2)}_{1'}\Ru_{1}
\nn\\ &=&
(\Ru_{2'}\tr \mu) \bigl(\Ru^{(1)}_{1'}\Ru_2\tr \la\bigr)\otimes\Ru^{(2)}_{1'}\Ru_{1} =
(\Ru_{2''}\Ru_{2'}\tr \mu) \bigl(\Ru_{1''}\Ru_2\tr \la\bigr)\otimes\Ru_{1'}\Ru_{1}
\nn\\ &=&
(\Ru_2\tr \la)(\Ru_{2'}\tr \mu)\otimes \Ru_{1'}\Ru_1.
\nn
\ee
In the last transformation, we have used the quasi-commutativity of
the algebra $\Rg$.

Employing the same arguments, we find
\be
 t (\mu) s (\la)&=&(\Ru_{2'}\tr \mu\otimes\Ru_{1'})(\la\otimes 1)= (\Ru_{2'}\tr \mu)
(\Ru^{(1)}_{1'}\tr \la)\otimes\Ru^{(2)}_{1'}
\nn\\
&=& (\Ru_{2''}\Ru_{2'}\tr \mu) (\Ru_{1''}\tr \la)\otimes\Ru_{1'}
=\la (\Ru_{2'}\tr \mu)\otimes \Ru_{1'}=
 s (\la) t (\mu).
\hspace{32pt}
\ee
We have proven the statement regarding the map $t^+$.
The case of $t^-$ is treated similarly, with $\Ru$ replaced by $\Ru^-$.
\end{proof}
\begin{propn}[\cite{Lu}]
\label{L-bialg}
The associative algebra $\Rg\rtimes \Ha$  is equipped with two
$\Rg$-bialgebroid structures, $\Rg_\pm\!\rtimes \Ha$, with the
source map $ s $ and the target map $t^\pm$ from Lemma \ref{lemma-s-t},
coproduct $\Delta(\la\otimes h):= (\la\otimes  h^{(1)})\tp_\Rg(1\otimes h^{(2)})$,
and the counit $\ve(\la\otimes h):=\la\:\ve(h)$, $\la\otimes h\in \Rg\rtimes \Ha$.
The anchor action of $\Rg\rtimes \Ha$ on $\Rg$ is given by
$(\mu\otimes h) \vdash\la = \mu (h\tr\la)$,
for $\la\in  \Rg$ and $\mu\otimes h\in \Rg\rtimes \Ha$.
\end{propn}

Remark that the bialgebroid structures $\Rg_\pm\!\rtimes \Ha$ on the
same associative algebra $\Rg\rtimes \Ha$ are determined solely by
the structures of the base algebra $\Lc_\pm$ on $\Lc$ (in other words, by the
$\Ha$-coactions).
By default, we understand by $\Lc \rtimes \Ha$ the bialgebroid $\Lc_+\! \rtimes \Ha$.

Each bialgebroid $\Rg_\pm\rtimes\Ha$ has a natural sub-bialgebroid.
To describe them, let us recall that the R-matrices $\Ru^\pm$ define two Hopf algebra maps from $\Ha^*_{op}$
to $\Ha$, by formulas (\ref{+-maps}).
Then $\Ru^\pm\in \Ha^{\sst\{\pm\}}\tp \Ha^{\sst\{\mp\}}\subset \Ha\tp \Ha$,
where $\Ha^{\sst\{\pm\}}$ are Hopf subalgebras in $\Ha$ that
are the images of $\Ru^\pm$.
Note that $\Rg_\pm$ is a base algebra over the Hopf algebra $\Ha^{\sst\{\pm\}}$, since
the coaction $\delta^\pm$ actually takes its values in $\Ha^{\sst\{\pm\}}\tp \Rg$,
see (\ref{d+-}).
The algebra $\Rg\rtimes\Ha$ contains $\Rg\rtimes \Ha^{\sst\{\pm\}}$ as subalgebras.
\begin{propn}
\label{+-bialg}
$\Rg\rtimes \Ha^{\sst\{\pm\}}$ are sub-bialgebroids in
$\Rg_\pm\rtimes \Ha$.
\end{propn}
\begin{proof}
The formula (\ref{t-d}) shows that the maps $t^\pm$
 take  values in $\Rg\rtimes\Ha^{\sst\{\pm\}}$.
Therefore $\Rg\rtimes\Ha^{\sst\{\pm\}}$ are $\Rg$-sub-bimodules
in $\Rg_\pm\rtimes\Ha$.
The coproduct in $\Rg_\pm\rtimes\Ha$  restricts to $\Rg\rtimes\Ha^{\sst\{\pm\}}$, thus
we conclude that $\Rg\rtimes\Ha^{\sst\{\pm\}}$ are sub-bialgebroids.
\end{proof}
\begin{remark}
\label{rm+-}
An arbitrary  Hopf algebra $\Ha$ is identified with
$(\D\Ha)^{\sst\{-\}}$, if the double $\D\Ha$ is equipped
with the quasitriangular structure $\Theta\in \Ha^*_{op}\tp \Ha\subset (\D\Ha)^{\tp2}$.
Given an $\Ha$-base algebra $\Lc$, one can build an $\Lc$-bialgebroid
$\Lc\rtimes\Ha\simeq \Lc\rtimes (\D\Ha)^{\sst\{-\}}\subset \Lc_-\!\rtimes \D\Ha$,
according to the Proposition \ref{+-bialg}.
The target map in $\Lc\rtimes\Ha$ is expressed through
the coaction by formula
\be
\label{t-d1}
t(\la)=\la^{[2]}\otimes \gm^{-1} (\la^{(1)})
,
\ee
as a specialization of (\ref{t-d}).
\end{remark}
\subsection{Quantum groupoid $\Ha_\Rg$}
In this subsection we build a quasitriangular bialgebroid $\Ha_\Rg$
as a quotient of $\Rg_\pm\rtimes \Ha$ by a certain biideal.
This quotient eliminates the distinctions between the two
bialgebroids $\Rg_+\rtimes \Ha$ and $\Rg_-\rtimes \Ha$.

To proceed with the study of the bialgebroid $\Rg_\pm\rtimes \Ha$, we need an algebraic construction
to be described next.
\begin{lemma}
Let $\phi$ be  an endomorphism of an associative algebra $\Bg$.
The left ideal $J_{\phi}$  generated by the image of the endomorphism $\phi-\id$
is a $\phi$-invariant two-sided ideal.
It is the minimal among $\phi$-invariant two-sided ideals such that the
endomorphism of $\Bg/J_\phi$ induced  by $\phi$ is identical.
\end{lemma}
\begin{proof}
The identity
$\bigl(\phi(a)-a\bigr) b=\bigl(\phi(ab)-ab\bigr)-\phi(a)\bigl(\phi(b)-b\bigr) $
being valid for any  pair $a,b\in \Bg$ shows that $J_\phi$ is a two-sided ideal.
The minimality property is obvious.
\end{proof}
\begin{lemma}
\label{J-phi}
a) The linear endomorphism $\phi\colon \Rg\rtimes \Ha\to \Rg\rtimes \Ha$
given by
\be
\label{phi}
\phi(\la\otimes h):= (\Ru_2\Ru_{1'})\tr \la\otimes \Ru_1\Ru_{2'}h.
\ee
is an algebra automorphism.
b) The ideal $J_\phi$ can be presented in the form $s(\Rg)(\phi-\id)(\Rg\rtimes \Ha)$.
c) As a two-sided ideal, $J_\phi$  is generated by the image of the map $t^+-t^-$, i.e. by the set
$(t^+-t^-)(\Rg)$.
\end{lemma}
\begin{proof}
Denote by $v$ the Drinfeld element $\Ru_1\gm(\Ru_2)\in \Ha$, \cite{Dr2}.
It satisfies the identities
\be
\Ru_{21}\Ru=\Delta(v^{-1})(v\tp v)
,\quad
v h  v^{-1}=\gm^{-2}(h),
\quad
h\in \Ha.
\label{Dr-el}
\ee
It is easy to check, using (\ref{Dr-el}), that the map $\phi_0\colon\la\otimes h\mapsto v\tr \la\otimes  v h v^{-1}$ is an
automorphism of the algebra $\Rg\rtimes \Ha$.
Then $\phi$ from (\ref{phi}) coincides with the composition of two automorphisms
$\mathrm{Ad}^{-1}(1\otimes v)\circ \phi_0$; this proves a).
Since $\phi$ is identical on $1\tp\Ha$, the image of $\phi-\id$ is invariant under
the left regular $\Ha$-action on $\Rg\rtimes \Ha$. Therefore $J_{\phi}$ can be presented
as a left $\Rg$-submodule generated by the image of $\phi-\id$; this proves b).
Remark that as a two-sided ideal, $J_\phi$ is generated by the image of
the map $(\phi-\id)\circ s$ or, in terms of the Drinfeld element, by the relations
\be
\label{v-rel}
(v\tr\la)\tp 1=(1\tp v)(\la\tp 1) (1\tp v)^{-1}, \quad \la \in \Lc.
\ee

Notice that $\phi\circ t^-=t^+$; this implies the equality
$t^- \equiv t^+ \mod J_\phi$ or, explicitly,
\be
\label{+-}
\Ru^+_2\la\otimes \Ru^+_1\equiv\Ru^-_2\la\otimes \Ru^-_1\mod J_\phi
,\quad \la \in \Rg.
\ee
On the other hand,
$
(\phi-\id)(\la\otimes 1)=\bigr(t^+-t^-\bigr)(\Ru^-_2\la)(1\otimes\Ru^-_1)
$,
hence $J_\phi$ lies in the ideal generated by $\bigr(t^+-t^-\bigl)(\Rg)$;
this proves c).
\end{proof}
Remark that the ideal $J_\phi$ is zero if $\Ha$ is triangular, i.e. $\Ru^+=\Ru^-$.
\begin{propn}
\label{q-triang}
The ideal $J_\phi$ is a biideal in both bialgebroids $\Lc_\pm\rtimes\Ha$. The quotient bialgebroid
$\Ha_\Rg:=\bigl(\Rg_+\rtimes \Ha\bigr)/J_\phi=\bigl(\Rg_-\rtimes \Ha\bigr)/J_\phi$
is quasitriangular, with the universal R-matrix
being the image of
\be
\label{R-matrix_eq}
(1\otimes \Ru_1)\tp_{\Rg_{op}} (1\otimes \Ru_2)
\ee
under the projection along $J_\phi$.
\end{propn}
\begin{proof}
By  Lemma \ref{J-phi}.b, $J_\phi$ is a left $\Rg$-module generated by the image of $\phi-\id$.
Applying the counit $\ve$ from Proposition \ref{L-bialg} to the formula (\ref{phi}),
we find  $\ve \circ (\phi-\id)=0$. Thus $J_\phi$ lies in the kernel of $\ve$.

Let us prove that $J_\phi$ is a biideal in $\Rg_+\rtimes \Ha$.
By Lemma \ref{J-phi}.c, $J_\phi$ is generated by the set $(t^-- t^+)(\Rg)$.
Therefore,  it is  sufficient to check
that $(\Delta\circ t^-)(\la)\equiv(\Delta\circ t^+)(\la)$, where
the symbol $\equiv$ means  equality modulo
$J_\phi\tp_\Rg \bigl(\Rg\rtimes \Ha\bigr)+\bigl(\Rg\rtimes \Ha\bigr)\tp_\Rg J_\phi$
for all $\la\in \Rg$. We have for $(\Delta\circ t^-)(\la)$
\be
\Delta(\Ru^-_2\tr \la\otimes \Ru^-_1)
\hspace{-1pt}&=&
\bigl(\Ru^-_{2} \Ru^-_{2'} \tr\la\otimes \Ru^-_{1}\bigr)
\tp_\Rg(1\otimes \Ru^-_{1'})
\equiv\bigl(\Ru^+_{2} \Ru^-_{2'} \tr\la\otimes \Ru^+_{1}\bigr)\tp_\Rg(1\otimes \Ru^-_{1'})
\nn\\&=&
 t ^+(\Ru^-_{2'}\tr \la)\tp_\Rg(1\otimes \Ru^-_{1'})=1\tp_\Rg s (\Ru^-_{2'}\tr \la)(1\otimes \Ru^-_{1'})
\nn\\&=&
1\tp_\Rg(\Ru^-_{2}\tr \la \otimes \Ru^-_{1})\equiv 1\tp_\Rg(\Ru^+_{2}\tr \la \otimes \Ru^+_{1}).
\nn
\ee
But the last expression is equal to $1\tp_\Rg t ^+(\la)=(\Delta\circ t ^+)(\la)$
since $\Delta$ is an $\Rg$-bimodule map. This proves that $J_\phi$ is a biideal
in $\Rg_+\rtimes \Ha$. This also implies that $J_\phi$ is a biideal
in $\Rg_-\rtimes \Ha$, in view of the symmetry $+\leftrightarrow -$.
The quotient $\bigl(\Rg_+\rtimes \Ha\bigr)/J_\phi$ is canonically
isomorphic to $\bigl(\Rg_-\rtimes \Ha\bigr)/J_\phi$ as a bialgebroid,
since $t^+\equiv t^- \mod J_\phi$.

Let us show that (\ref{R-matrix_eq}) is a universal $R$-matrix in $\Ha_{\Rg}$.
Only condition 1 of Definition \ref{R-matrix} requires verification.
The other conditions follow from the properties of $\Ru$ as
a universal R-matrix of the Hopf algebra $\Ha$.

Computing the element $(1\otimes \Ru_1) t^+(\la)\tp_{\Rg_{op}} (1\otimes \Ru_2)$  modulo the ideal $J_\phi$
we find
\be
(1\otimes \Ru_1)(\Ru_2 \tr\la\otimes \Ru_1)\tp_{\Rg_{op}} (1\otimes \Ru_2)
&=&
(\Ru^{(1)}_1\Ru_{2'} \tr\la\otimes \Ru^{(2)}_1\Ru_{1'})\tp_{\Rg_{op}} (1\otimes \Ru_2)
\nn\\
&=&
 s (\Ru^{(1)}_1\Ru_{2'} \tr\la)(1\otimes \Ru^{(2)}_1\Ru_{1'})\tp_{\Rg_{op}} (1\otimes \Ru_2)
\nn\\
&=&
(1\otimes \Ru^{(2)}_1\Ru_{1'})\tp_{\Rg_{op}}  t ^+(\Ru^{(1)}_1\Ru_{2'} \tr\la)(1\otimes \Ru_2).
\nn
\ee
Since $t^+\equiv t^- \mod J_{\phi}$, this expression becomes equal to
\be
(1\otimes \Ru_{1''}\Ru_{1'})\tp_{\Rg_{op}} (\Ru_{2'}\tr \la\otimes \Ru_{2''})
&=&
(1\otimes \Ru_{1'}\Ru_1)\tp_{\Rg_{op}} (\Ru_2\tr \la\otimes \Ru_{2'})
\nn
\ee
and, finally, to $(1\otimes \Ru_1)\tp_{\Rg_{op}} (1\otimes \Ru_2) s (\la)$, as required.
Let us comment that one can  also deduce this fact directly from Remark \ref{R-tw},
noticing  that the inverse R-matrix of $\Ha$ is a twisting cocycle of $\Bg$.
\end{proof}
\noindent
\begin{remark}
In what follows we will abuse notation suppressing the projection $\Rg_\pm\rtimes \Ha\to \Ha_\Rg$ when
writing elements of $\Ha_\Rg$. In other words, the reader can perceive calculations
in $\Ha_\Rg$ as those in $\Rg_\pm\rtimes \Ha$ done modulo $J_\phi$. The most important
feature for us is the identity $\Ru^+_2\la\otimes \Ru^+_1=\Ru^-_2\la\otimes \Ru^-_1$,
which is valid in $\Ha_\Rg$ for all $\la \in \Rg$.
\end{remark}
\section{On the antipode}
\label{secAntipode}
In this subsection we study antipodes in bialgebroids $\Lc\rtimes\Ha$ and $\Ha_\Lc$.
It turns out that they can be defined as isomorphisms between opposite and coopposite
bialgebroids, analogously to Hopf algebras.
However, contrary to the Hopf algebra case, there is no canonical way
to define the opposite bialgebroid. Even the coopposite bialgebroid,
although defined canonically in  Subsection \ref{subsecQST}, is in fact over the opposite base.
Nevertheless, using the specific form of the bialgebroids under consideration,
the opposite bialgebroids can be introduced.

\subsection{Bialgebroid $(\Lc\rtimes\Ha)_{op}$}
\label{subsec_ba_op}
In this and the next subsections we consider the $\Lc$-bialgebroid
$\Lc\rtimes\Ha$, where
$\Lc$ is a base algebra for a general (not necessarily quasitriangular)
Hopf algebra  $\Ha$, cf. Remark \ref{rm+-}.

\begin{lemma}
\label{op}
Let $\Ha$ be a Hopf algebra. Let $\Lc$ be an $\Ha$-base algebra with the $\Ha$-action $\tr$ and the
coaction $\delta$. Then $\Lc_{op}$ is a base algebra over $\Ha_{op}$ with respect to the $\Ha_{op}$-action
$x\btr\ell:= \gm^{-1}(x)\tr \ell$ and the $\Ha_{op}$-coaction $\delta^\op=\delta$.
\end{lemma}
\begin{proof}
Obviously  $\Lc_{op}$ is a left $\Ha_{op}$-module algebra with respect to $\btr$ and a left $\Ha_{op}$-comodule algebra
with respect to $\delta^\op$.
Let us show that it is a Yetter-Drinfeld module  with respect to $\Ha_{op}$.
This is equivalent to the condition
$$
\delta(x\btr \ell)=x^{(1)}\cdot\ell^{(1)}\cdot\gm^{-1}(x^{(3)})\tp  x^{(2)}\btr \ell^{[2]},
$$
which is the formula (\ref{ba1}) translated to the case of $\Lc_{op}$ and $\Ha_{op}$ instead
of $\Lc$ and $\Ha$ (note
that $\gamma^{-1}$ the antipode for $\Ha_{op}$). The dots mean the opposite multiplication.

Finally, the $\Ha$-commutativity condition (\ref{ba2}) in $\Lc$ transforms into
the $\Ha_{op}$-commutativity condition $\la\cdot\mu=(\la^{(1)}\btr\mu)\cdot\la^{[2]}$ in $\Lc_{op}$.
\end{proof}

Consider the opposite associative algebra $(\Lc\rtimes\Ha)_{op}$. This algebra
contains $\Lc_{op}$ and $\Ha_{op}$ as subalgebras and, in fact, has
the form of smash product $\Lc_{op}\rtimes\Ha_{op}$, where the action of $\Ha_{op}$
on $\Lc_{op}$ is specified in Lemma \ref{op}. But
$\Lc_{op}$ is a base algebra over $\Ha_{op}$, by Lemma \ref{op}, thus
$\Lc_{op}\rtimes\Ha_{op}$ is equipped by the structure of  $\Lc_{op}$-bialgebroid
in the standard way.
Thus there is a canonical bialgebroid structure on the
opposite algebra $(\Lc\rtimes\Ha)_{op}$.
We denote by  $s^\op$, $t^\op$, $\Delta^\op$, and
$\ve^\op$ respectively, the source, target,
coproduct, and counit maps of the opposite bialgebroid $\Lc_{op}\rtimes\Ha_{op}$.
\begin{definition}
\label{smash_op}
The opposite bialgebroid $(\Lc\rtimes \Ha)_{op}$ to
the bialgebroid  $\Lc\rtimes \Ha$
is an $\Lc_{op}$-bialgebroid $(\Lc_{op}\rtimes\Ha_{op},\Lc_{op},s^\op,t^\op,\Delta^\op,\ve^\op)$
\end{definition}

\subsection{The antipode in $\Lc\rtimes\Ha$}
Recall that  $\Theta\in \Ha^*_{op}\tp\Ha
\subset (\D\Ha)^{\tp 2}$ denotes the standard quasitriangular structure on $\D\Ha$.
Let us compute the target map of the bialgebroid $(\Lc\rtimes \Ha)_{op}$
defined in the previous subsection, in terms of $\Lc$, $\Ha$, and $\Theta$.
\begin{lemma}
\label{target_op}
The map $\la\mapsto \Theta_1\tr\la\tp \gm(\Theta_2)$ from
$\Lc$ to $\Lc\tp \Ha$ yields the target map $t^\op$ of the $\Lc_{op}$-bialgebroid $\Lc_{op}\rtimes \Ha_{op}$.
\end{lemma}
\begin{proof}
Let us apply the formula (\ref{t-d1})
to the bialgebroid $\Lc_{op}\rtimes \Ha_{op}$.
The coaction $\delta^\op$ coincides with $\delta$, which
has the form $\la\mapsto \Theta_2\tp \Theta_1\tr \la$.
Now the lemma follows from the fact that the antipode in $\Ha_{op}$ is
the inverse antipode in $\Ha$.
\end{proof}

Consider the map $\zeta\colon\Lc\rtimes\Ha\to \Lc_{op}\rtimes\Ha_{op}$ defined by
\be
\label{op-coop}
\zeta\colon\la\otimes h\mapsto \Theta_1\tr\la\otimes \gm(\Theta_2)\cdot\gm(h),
\ee
where $\gm$ is the antipode of $\Ha$ and $\Theta$ is the R-matrix of $\D\Ha$
(here we suppressed the anti-isomorphism $\iota\colon\Lc\to \Lc_{op}$).
\begin{propn}
The map (\ref{op-coop}) defines an isomorphism of $\Lc_{op}$-bialgebroids
\be
\label{inv_ap}
\zeta\colon(\Lc\rtimes\Ha)^{op}\to \Lc_{op}\rtimes\Ha_{op}.
\ee
\end{propn}
\begin{proof}
The map (\ref{op-coop}) is an algebra homomorphism
when  restricted to the subalgebras $\Lc\tp 1$ and $1\tp \Ha$.
By construction, it respects the product $(\la\otimes 1)(1\otimes h)$.
To complete the proof,
one must check that $\zeta$ respects the product $(1\otimes h)(\la\otimes 1)$:
\be
\zeta(1\otimes h)\zeta(\la\otimes 1)&=&
\bigl(1\otimes \gm(h)\bigr)(\Theta_1\tr \la\otimes \Theta_2)=\gm(h^{(2)})\btr\Theta_1\tr \la\otimes \gm(h^{(1)})\cdot\gm(\Theta_2)
\nn\\
&=&
h^{(2)}\Theta_1\tr \la\otimes \gm(h^{(1)})\cdot\gm(\Theta_2)=
h^{(2)}\Theta_1\tr \la\otimes \gm(h^{(1)}\Theta_2).
\nn
\ee
Since $\Theta$ is an R-matrix of the double $\D\Ha$, the above expression can be rewritten
as
$$
(\Theta_1 h^{(1)}\tr)\la\otimes \gm(\Theta_2)\cdot\gm(h^{(2)})=\zeta(h^{(1)}\tr\la\otimes h^{(2)})
=\zeta\bigl((1\otimes h)(\la\otimes 1)\bigr).$$

The target map  $t^{op}$ of the  $\Lc_{op}$-bialgebroid $(\Lc\rtimes\Ha)^{op}$
comes from the source map of the $\Lc$-bialgebroid $\Lc\rtimes\Ha$ and it is
equal to $s\circ\iota^{-1}$.
It follows from Lemma \ref{target_op} that $\zeta\circ t^{op}=t^\op$.
Let us prove that $\zeta\circ s^{op}$ equals
the source map $s^\op$ of $\Lc_{op}\rtimes\Ha_{op}$. Suppressing the notation of the map $\iota$,
we have
\be
(\zeta\circ s^{op})(\la)
&=&
(\zeta\circ t)(\la)=\zeta\bigl(\Theta_1\tr\la\otimes \gm^{-1}(\Theta_2)\bigr)
=
(\Theta_{1'}\Theta_1)\tr\la\otimes \gm(\Theta_{2'})\cdot\Theta_2 =\la\otimes 1.
\nn
\ee
Here we have used the fact
$\Theta_{1'}\Theta_1\tp \Theta_2\gm(\Theta_{2'})=1\tp 1$, which is
equivalent to the standard identity
$\Theta_1\tp \gm^{-1}(\Theta_2)=\Theta^{-1}$ for the universal R-matrix.
Thus we have shown that $\zeta$ is a morphism of $\Lc_{op}$-bimodules.

Now let us show that $\zeta$ respects the coproducts.
Indeed, we have
$$
\bigl((\zeta\tp_{\Lc_{op}}\zeta)\circ \Delta^{op}\bigr)(\la\otimes h )
=
\zeta(1\otimes h^{(2)})\tp_{\Lc_{op}} \zeta(\la\otimes h^{(1)})=\bigl(1\otimes \gm(h^{(2)})\bigr)\tp_{\Lc_{op}} t^\op(\la)\gm(h^{(1)}).
$$
The rightmost expression is equal to $\Delta^\op\bigl(t^\op(\la)\gm(h)\bigr)=\Delta^\op\bigl(\zeta(\la\otimes h)\bigr)$.
Thus we have checked the
right equation from (\ref{bial_hom}). The left one, concerning the counits, readily
follows from the definition of $\zeta$.
\end{proof}

Replacing $\Lc$ and $\Ha$ by $\Lc_{op}$ and $\Ha_{op}$  and taking the inverse map in (\ref{inv_ap}),
we obtain a bialgebroid isomorphism
\be
\label{antipode}
\Lc\rtimes\Ha\to (\Lc_{op}\rtimes\Ha_{op})^{op}.
\ee
Using the argument after the proof of Lemma \ref{op}, we can consider the map
(\ref{antipode}) as an anti-isomorphism of the
associative algebra $\Lc\rtimes\Ha$, which reads
\be
\label{antipode_exp}
\la\otimes h\mapsto \bigl(1\otimes\gamma(h)\bigr)t( v^{-1}\tr \la),\quad \la\tp h\in \Lc\rtimes\Ha,
\ee
where $ v$ is the Drinfeld element, cf. (\ref{Dr-el}).
We denote the map  (\ref{antipode}) by $\gm$ regarding it
as an extension of the antipode of $\Ha\subset \Lc\rtimes \Ha$.
\begin{remark}
The map (\ref{antipode_exp}) coincides with the antipode of  \cite{Lu}.
However it was considered there just as an operator on $\Lc\rtimes\Ha$
possessing a certain set properties. We would like
to emphasize the bialgebroid meaning of the map
(\ref{antipode_exp}). Namely,
it implements  an isomorphism (\ref{antipode})  between {\em different}
bialgebroids over  {\em different} bases. This gives rise to the categorical interpretation
of the antipode  of \cite{Lu}. It gives rise to  an isomorphism of the corresponding
monoidal categories of modules.
\end{remark}

\subsection{The antipode in the quantum groupoid $\Ha_\Rg$}
In this subsection we will investigate the behavior of the antipode (\ref{antipode}) under
the projection $\Rg\rtimes \Ha\to \Ha_\Rg$  assuming $\Ha$ quasitriangular with R-matrix $\Ru$
and $\Rg$  quasi-commutative; $\Rg$ is equipped with the $\Ha$-base algebra structure $\Rg_+$, cf. Remark \ref{rm+-}.
Denote by $v^\op$ the Drinfeld element of $\Ha_{op}$
and by $\phi^\op$ the automorphism (\ref{phi})
specialized to the case of the bialgebroid $\Rg_{op}\rtimes\Ha_{op}$.
According to Proposition \ref{q-triang}, the ideal $J_{\phi^\op}$
is a biideal in $\Rg_{op}\rtimes \Ha_{op}$.

Thus we have two two-sided ideals in the algebra $\Rg\rtimes\Ha$,
namely $J_{\phi}$ and $J_{\phi^\op}$, corresponding to the
bialgebroids $\Rg\rtimes \Ha$ and $\Rg_{op}\rtimes \Ha_{op}$
(two-sided ideals are the same in opposite algebras).
\begin{lemma}
\label{lem_ideal_op}
The ideals $J_{\phi}$ and $J_{\phi^\op}$
coincide.
\end{lemma}
\begin{proof}
In the course of the proof of Proposition \ref{q-triang} we have shown that
the ideal $J_\phi$ is generated by the relations
(\ref{v-rel}). Specializing (\ref{v-rel}) for the bialgebroid
$\Rg_{op}\rtimes\Ha_{op}$, we find that the ideal $J_{\phi^\op}$ is
generated by
the relations
\be
\label{phi_bul}
v^\op\btr \la\tp 1=(1\tp v^\op)\cdot(\la\tp 1)\cdot(1\tp v^\op)^{-1},
\ee
where the dots stand for the multiplication in $\Rg_{op}\rtimes\Ha_{op}$.

Observe that the antipode in $\Ha_{op}$ is the inverse
antipode in $\Ha$.
Therefore $v^\op=v^{-1}$ and $\gm^{-1}(v^{-1})=v$, since  $v^{-1}$ implements
the squared antipode by conjugation.  Taking into account this argument,
equation (\ref{phi_bul}) translates into equation (\ref{v-rel}) in $\Rg\rtimes\Ha$.
\end{proof}

Following Proposition \ref{q-triang} we can introduce a  quantum groupoid that
is opposite to $\Ha_\Lc$.
\begin{definition}
\label{op_QGd}
The opposite quantum groupoid $(\Ha_\Rg)_{op}$ is
a quasitriangular  $\Rg_{op}$-bialgebroid that is the
quotient of
$\Rg_{op}\rtimes \Ha_{op}$ by the biideal $J_{\phi^\op}$.
\end{definition}
The coopposite $\Rg$-bialgebroid $(\Ha_\Rg)^{op}_{op}$
to the $\Rg_{op}$-bialgebroid $(\Ha_\Rg)_{op}$ is defined canonically,
see Subsection \ref{subsecQST}.
\begin{propn}
The antipode (\ref{antipode}) descends to an isomorphism of quantum groupoids
$\gamma:\Ha_{\Rg}\to (\Ha_{\Rg})^{op}_{op}$.
\end{propn}
\begin{proof}
By Lemma \ref{lem_ideal_op}, the ideal $J_\phi$ defining
$\Ha_{\Rg}$ coincides with the ideal defining $(\Ha_{\Rg})_{op}$
and thus $(\Ha_{\Rg})^{op}_{op}$. Therefore, it suffices to
check that $J_\phi$ is invariant with respect to the antipode (\ref{antipode}).
Identifying $\Rg$ and $\Ha$ with the corresponding subalgebras in $\Rg\rtimes \Ha$,
we can write
\be
\gm(v\la v^{-1})&=&
v\gm(\la)v^{-1}=
v t^+(v^{-1}\tr\la)v^{-1}
=
\bigl((v^{(1)}\Ru_2 v^{-1})\tr\la\bigr)(v^{(2)}\Ru_1 v^{-1}).
\label{eq_2}
\ee
for any $\la\in \Rg$.
Using (\ref{Dr-el}) and the identity $(\gm\tp \gm)(\Ru)=\Ru$, we
find the last expression to be equal to $(\Ru^{-1}_1\tr\la)\Ru_2^{-1}=t^-(\la)\equiv t^+(\la) \mod J_\phi$
and therefore to $\gm(v\tr\la)$ modulo the ideal $J_\phi$,
by formula (\ref{antipode_exp}) for $h=1$.
Thus the relations (\ref{v-rel}) are preserved by $\gm$, modulo $J_{\phi}$.

The induced homomorphism $\gamma \colon\Ha_{\Rg}\to(\Ha_{\Rg})^{op}_{op}$ of bialgebroids
relates the quasitriangular structures
of $\Ha_{\Rg}$ and $(\Ha_{\Rg})^{op}_{op}$, i.e. $\gamma\tp\gamma$ leaves the
R-matrix invariant.
Thus $\gamma$ is an isomorphism of quantum groupoids.
\end{proof}
\section{Dynamical cocycles and twisting bialgebroids}
\label{secBTDC}
\subsection{Twisting by dynamical cocycles}
\label{subsecTDC}
The present section establishes a relation between bialgebroid twists
and dynamical cocycles over a non-abelian base from \cite{DM1}.
The case of abelian base was treated in \cite{Xu1}.

A categorical definition of dynamical twist is given in \cite{DM1}.
Here we will work with the equivalent definition in terms of universal dynamical twisting
cocycle.
A \select{universal dynamical cocycle} over an $\Ha$-base algebra $\Lc$
is
an invertible element $\F=\F_1\tp\F_2\tp\F_3\in \U\tp \U\tp\Lc $, where $\U$ is a Hopf algebra containing $\Ha$,
satisfying the invariance condition
\be\label{coninv1}
 h^{(1)}\F_1\tp h^{(2)}\F_2\tp h^{(3)}\tr \F_3
 &=&
 \F_1 h^{(1)}\tp \F_2 h^{(2)}\tp \F_3,\quad \forall h\in \Ha,
\ee
the shifted cocycle condition
\be
\label{shifted_cocycle}
(\Delta\tp \id)(\F)\;\bigl(\F_1\tp\F_2\tp \F^{(1)}_3\tp \F^{[2]}_3\bigr)
&=&
(\id\tp \Delta)(\F)(\F_{23}),
\ee
and the normalization condition
\be
(\ve \tp \id\tp \id)(\F)=&1\tp 1\tp 1& = (\id\tp \ve \tp \id)(\F).
\ee
Note that equation (\ref{shifted_cocycle}) holds in $\U\tp\U\tp \U\tp \Lc$.

Assume that $\Ha$ is quasitriangular with the R-matrix $\Ru$ and
$\Lc$ is $\Ha$-commutative. Recall that $\Lc$ can be equipped
with two $\Ha$-base algebra structures $\Lc_\pm$ by the coactions
(\ref{d+-}).
Consider the tensor product bialgebroid $\U\tp(\Lc_+\rtimes \Ha)$,
as in Example \ref{H-B}.
\begin{propn}
\label{dc-tw}
Let $\F=\F_1\tp \F_2\tp \F_3\in \U\tp \U \tp \Lc_-$ be a dynamical twist.
Then the element $\Psi\in (\U\tp \Lc_+\rtimes\Ha)\tp_\Lc (\U\tp \Lc_+\rtimes\Ha)$,
\be
\label{biatwist}
\Psi&:=&(\F_1\tp \F_{3} \otimes \Ru_1)\tp_{\Lc}(\F_2\Ru_2\tp 1\otimes 1),
\ee
 is a bialgebroid  twist.
\end{propn}
\begin{proof}
Explicitly, the comultiplication in $\U\tp (\Lc_+\rtimes\Ha)$ is
written as
\be
\Delta(u\tp \la\otimes h)&:=&(u^{(1)}\tp \la\otimes  h^{(1)})\tp_{\Lc}
(u^{(2)}\tp 1\otimes h^{(2)})
\ee
for any $u\tp \la \otimes h\in \U\tp (\Lc\rtimes\Ha)$.
Then the right-hand side of (\ref{bialg_cocycle})
 is equal to
\be
\label{eq1}
(\F_1\tp \F_{3}\otimes \Ru_1)\tp_{\Lc} (\F^{(1)}_2\Ru^{(1)}_2 \F_{1'}\tp \F_{3'} \otimes \Ru_{1'''})\tp_{\Lc}
(\F^{(2)}_2\Ru^{(2)}_2\F_{2'}\Ru_{2'''}\tp 1 \otimes 1).
\ee
By the standard Hopf algebra technique, the identity
(\ref{coninv1}) implies
$$
h^{(1)}\F_1\tp h^{(2)}\F_2\otimes \F_3=\F_1 h^{(1)}\tp \F_2 h^{(2)}\otimes \gm^{-1}(h^{(3)})\tr\F_3,
\quad \forall h\in \Ha,
$$
where  $\gm$ is the antipode in $\Ha$.
Using this, we transform (\ref{eq1}) to
\be
(\F_1\tp \F_{3} \otimes \Ru_1)\tp_{\Lc} (\F^{(1)}_2\F_{1'} \Ru^{(1)}_2 \tp \gm^{-1}\Ru^{(3)}_2\tr\F_{3'}\otimes
\Ru_{1'''})\tp_{\Lc}
(\F^{(2)}_2\F_{2'}\Ru^{(2)}_2\Ru_{2'''}\tp 1 \otimes 1) =\nn
\hspace{-6pt}\\
(\F_1\tp \F_{3} \otimes \bar\Ru_{1''}\Ru_1)\tp_{\Lc}
(\F^{(1)}_2\F_{1'} \Ru^{(1)}_2 \tp \bar\Ru_{2''}\tr\F_{3'}\otimes
\Ru_{1'''})\tp_{\Lc}
( \F^{(2)}_2\F_{2'}\Ru^{(2)}_2\Ru_{2'''}\tp 1 \otimes 1).
\nn
\ee
The element $\bar \Ru$ denotes
$\Ru_{1}\tp \gm^{-1}\Ru_{2}$, which is the inverse to $\Ru$.
The term $\bar\Ru_{2''}\tr\F_{3'}$ in the middle tensor factor can
be pulled to the left as the factor $s(\bar\Ru_{2''}\tr\F_{3'})$.
Using the definition of tensor product  over $\Lc$
we transform this expression  to
\be
 t^+ (\bar\Ru_{2''}\tr\F_{3'})\Bigl(\F_1\tp\F_{3} \otimes\bar\Ru_{1''}\Ru_1\Bigr)
\tp_{\Lc}
\Bigl(\F^{(1)}_2\F_{1'}\Ru^{(1)}_2 \tp 1\otimes   \Ru_{1'''}\Bigr)\tp_{\Lc} \ldots
\nn \hspace{85pt}\\
=\Bigl(1\tp (\Ru_{2''}\bar\Ru_{2})\tr\F_{3'}\otimes \Ru_{1''}\Bigr)
\Bigl(\F_1\tp \F_{3} \otimes \bar\Ru_{1}\Ru_1\Bigr)\tp_{\Lc}
\Bigl(\F^{(1)}_2\F_{1'} \Ru^{(1)}_2 \tp 1\otimes \Ru_{1'''}\Bigr)\tp_{\Lc} \ldots
\nn \\
=
\Bigl(\F_1\tp \F_{3}\bigl((\Ru_{2''}\bar\Ru_{2})\tr\F_{3'}\bigr) \tp \Ru_{1''}\bar\Ru_{1}\Ru_1\Bigr)\tp_{\Lc}
\Bigl(\F^{(1)}_2\F_{1'} \Ru^{(1)}_2 \tp 1\otimes \Ru_{1'''}\Bigr)\tp_{\Lc} \ldots
\nn \hspace{38.pt}\\
=(\F_1\tp \F_{3}\F_{3'}\otimes\Ru_{1''}\Ru_{1'} )
\tp_{\Lc}
(\F^{(1)}_2\F_{1'} \Ru_{2'} \tp 1\otimes \Ru_{1'''} )\tp_{\Lc}
(\F^{(2)}_2\F_{2'}\Ru_{2''}\Ru_{2'''}\tp 1 \otimes 1).
\nn
\hspace{13.5pt}
\ee
Here we employed the fact that the image of the map $t^+$ commutes with all the elements $(x\tp \mu \tp 1)\in \U\tp \Lc_+\rtimes\Ha$.

On the other hand, the left-hand side of (\ref{bialg_cocycle}) turns to
\be
\Bigl(\F^{(1)}_1\F_{1'}\tp \F_{3}(\Ru_{1}\tr\F_{3'})\otimes \Ru_{1''}\Ru_{1'}\Bigr)
\tp_{\Lc}\hspace{8cm}
\nn\\
\tp_{\Lc}
\Bigl(\F^{(2)}_1\F_{2'}\Ru_{2'}\tp 1\otimes \Ru_{1'''}\Bigr)\tp_{\Lc}
\Bigl(\F_2\Ru_{2}\Ru_{2''}\Ru_{2'''}\tp 1 \otimes  1\Bigr).
\nn
\ee
Thus $\Psi$ satisfies equation (\ref{bialg_cocycle}) if $\F$ is a dynamical cocycle
over the base algebra $\Lc_-$, with the coaction
$\delta^-(\la)=\Ru_2\tp \Ru_1\tr \la$.
\end{proof}
\begin{corollary}
A dynamical cocycle $\F\in \U\tp \U\tp \Lc_-$ defines a new $\Lc$-bialgebroid
structure  $\widetilde{\U\tp (\Lc_+\rtimes\Ha)}$ on
the algebra $\U\tp (\Lc_+\rtimes\Ha)$,
with the same counit and  target map
$\tilde  t := t $
but the new source map
$\tilde  s (\la):= \Ru_2\tp \Ru_1 \tr\la\otimes 1$
and comultiplication (\ref{twcopr}) with $\Psi$ given by
(\ref{biatwist})
\end{corollary}
\begin{proof}
This is a corollary of Proposition \ref{twist}
and Proposition \ref{dc-tw}.
The source and target maps are readily calculated.
\end{proof}

In fact, the twist $\Psi$ is supported in the sub-bialgebroid $\U\tp (\Lc\rtimes \Ha^{\sst\{+\}})$.
Therefore the sequence  of bialgebroid homomorphisms
$$
\U\tp \Lc\rtimes\Ha^{\sst\{+\}}\to\U\tp \Lc_+\rtimes\Ha\to
\U\tp \Ha_\Lc,
$$
where the left arrow is embedding and the right one is projection,
gives rise to the sequence of homomorphisms of twisted bialgebroids, cf. Remark
\ref{hom_tw}:
$$
\widetilde{\U\tp \Lc\rtimes\Ha^{\sst\{+\}}}\to\widetilde{\U\tp \Lc_+\rtimes\Ha}
\to
\widetilde{\U\tp \Ha_\Lc}.
$$
\begin{corollary}
\label{corUDHL}
Suppose $\U$ is quasitriangular with the universal R-matrix
$\Omega$. Then the bialgebroid
$\widetilde{\U\tp \Ha_\Lc}$ is quasitriangular, with the universal R-matrix
\be
\label{R-matrix-twisted}
(\Ru^{-1}_{2'}\Ru^{-1}_{2'''}\tilde \Omega_1\tp 1\otimes \Ru_1\Ru_{1''})
\tp_{\Lc_{op}}
(\tilde \Omega_2\Ru_{2''}\tp \Ru^{-1}_{1'''}\tr\tilde \Omega_3\otimes\Ru^{-1}_{1'}\Ru_2),
\ee
where $\tilde \Omega:= \F_{21}^{-1} \Omega \F=\tilde \Omega_1\tp \tilde \Omega_2\tp \tilde \Omega_3\in \U\tp \U\tp \Lc$.
\end{corollary}
\begin{proof}
If $\U$ is quasitriangular with the universal R-matrix $\Omega$, then the tensor product bialgebroid $\U\tp \Ha_\Lc$
is also quasitriangular with the universal R-matrix
\be
(\Omega_1\tp 1\otimes\Ru_1)\tp_{\Lc_{op}} (\Omega_2\tp 1\otimes \Ru_2)
.
\label{R-matrix-triv}
\ee
The quasitriangular structure on $\widetilde{\U\tp \Ha_\Lc}$
is obtained from (\ref{R-matrix-triv})
by twisting, following Proposition \ref{Rtwisted}.
It can be expressed through the element $\tilde \Omega$ by the formula
(\ref{R-matrix-twisted}).
Let us remark that the element $\tilde \Omega$ is a solution to
the DYBE over the base algebra $\Lc$, see \cite{DM1}.
\end{proof}
\subsection{The twisted tensor product \tw{\U}{\Ru}{\Ha_\Lc}}
\label{subsecTTP}
The element $\F=1\tp 1\tp 1$ is a particular case of dynamical cocycle.
So we can always build a twist by (\ref{biatwist}) $\Psi_\Ru:=\Psi|_{\F=1}$.
A slight modification of the proof of Proposition \ref{dc-tw}
shows that the bialgebroid $\U\tp(\Lc_+\rtimes \Ha)$ has one more
twist, namely if $\Ru$ in $\Psi_\Ru$ is replaced by $\Ru^-$.
The same holds true for the bialgebroid $\U\tp(\Lc_-\rtimes \Ha)$,
which differs from $\U\tp(\Lc_+\rtimes \Ha)$ by the alternative choice
of the quasitriangular structure on $\Ha$.
These twists are analogous to the twisted tensor
products of Hopf algebras, cf. Subsection \ref{subsecHAD} and Example \ref{H-B}.
Following this analogy,
we reserve the special notation,
\tw{\U}{\Ru}{(\Lc\rtimes\Ha)} and \tw{\U}{\Ru}{\Ha_\Lc} of the bialgebroids
$\U\tp (\Lc\rtimes\Ha)$ and $\U\tp \Ha_\Lc$ twisted with $\Psi_\Ru$.

The goal of the present subsection is to establish the following commutative diagram of
bialgebroid homomorphisms:
\be
\begin{array}{ccccc}
\Lc_-\rtimes \Ha&\longrightarrow &\mbox{\tw{\Ha}{\;\;\Ru^+}{(\Lc_-\rtimes\Ha)}}
&\longrightarrow &\mbox{\tw{\U}{\;\;\Ru^+}{(\Lc_-\rtimes\Ha)}}\\
\downarrow&&\downarrow&&\downarrow\\
\Ha_\Lc& &\mbox{\tw{\Ha}{\;\;\Ru^+}{\Ha_\Lc}}
&\longrightarrow &\mbox{\tw{\U}{\;\;\Ru^+}{\Ha_\Lc}}\\
\end{array}
\ee
The similar diagram takes place upon interchanging $+\leftrightarrow -$.
The horizontal arrows on the right are obvious: the twist of the rightmost
bialgebroids is transferred from the middle ones via bialgebroid
homomorphisms, cf. Remark \ref{hom_tw}.
The horizontal arrow on the left  descends from
the coproduct $\Ha\to \Ha\tp \Ha$. Thus it can be viewed as
a generalization of the Hopf algebra embedding
$\Ha\stackrel{\Delta}{\longrightarrow} \mbox{\tw{\U}{\;\;\Ru^\pm}{\Ha}}$.
The blank space on the left of the bottom line means that
there is no homomorphisms from $\Ha_\Lc$ to $\mbox{\tw{\Ha}{\;\;\Ru^+}{\Ha_\Lc}}$
in general.
One can construct  a homomorphism from $\Ha_\Lc$ into the quotient
of  $\mbox{\tw{\Ha}{\;\;\Ru^+}{\Ha_\Lc}}$ by the ideal
generated by the image
of $J_\phi\subset \Lc_-\rtimes \Ha$.
Note that the two-sided ideal generated by the image of
a biideal under a bialgebroid homomorphism is always a biideal.
We do not focus on this issue here.

\begin{propn}
\label{l_hom}
Let $\Ha$ be a quasitriangular Hopf algebra with the R-matrix
$\Ru$ and $\Lc$ is an $\Ha$-commutative algebra
considered as a base algebra $\Lc_-$.
The map
\be
\label{hom}
\eta\colon \la\otimes h\mapsto \Ru_2 h^{(1)}\tp \Ru_1\tr\la \otimes h^{(2)}
\ee
from $\Lc_-\rtimes\Ha$ to \tw{\Ha}{\Ru^+}{(\Lc_-\rtimes\Ha)} is a bialgebroid embedding.
\end{propn}
\begin{proof}
Let us prove that (\ref{hom}) is an algebra homomorphism.
When restricted to $\Lc$, $\eta$ coincides with the source map $\tilde s$, whereas
the restriction to $\Ha$ descends from the coproduct of $\Ha$. Thus $\eta$ is a homomorphism
on the subalgebras $\Lc$ and $\Ha$ in $\Lc\rtimes \Ha$.
We have  $\eta(\la\tp 1)\eta(1\tp h)=\eta\bigl((\la\tp 1)(1\tp  h)\bigr)=\eta(\la\tp h)$
by construction. Further,
$$
\eta\bigl((1\tp  h)(\la\tp 1)\bigr)=\Ru_2 h^{(2)}\tp \Ru_1 h^{(1)}\tr\la\otimes h^{(3)}=
 h^{(1)}\Ru_2\tp h^{(2)}\Ru_1 \tr\la \otimes h^{(3)}=\eta(1\tp h)\eta(\la\tp 1).
$$
This proves that $\eta$ is an algebra homomorphism.
It is an embedding, since there is a projection
$\ve_\Ha\tp \id\colon\mbox{\tw{\Ha}{\Ru}{(\Lc_-\rtimes\Ha)}}\to \Lc_-\rtimes\Ha$
(in fact, this is a bialgebroid map) and
the composition of $\eta$ with this projection is identical on $\Lc_-\rtimes\Ha$.

Let us show that $\eta$ is an $\Lc$-bimodule map.
The equality $\tilde s= \eta\circ s$ is obvious,
so let us consider
the target maps. We have, for $\la\in \Lc$,
$$
(\eta\circ t^-)(\la)=
\eta(\Ru^-_2\tr\la\otimes \Ru^-_1)=
 \Ru_{2'}\Ru^-_{1''}\tp\Ru_{1'}\Ru^-_{2''}\Ru^-_{2}\tr\la\otimes\Ru^-_{1}=
1\tp\Ru^-_{2}\tr\la\otimes\Ru^-_{1}
=\tilde t^-(\la).
$$
Here we used $\Ru^-=\Ru_{21}^{-1}$.

We must show that $\eta$ respects the coproducts.
It is obvious for $\eta$ restricted to $\Lc\subset\Lc_- \rtimes \Ha$, so it suffices to check this on the elements of
$\Ha\subset \Lc_- \rtimes \Ha$.
When restricted to $\Ha$, the map $\eta$
coincides with
$\Delta\colon\Ha\to \Ha \tp \Ha$.
The bialgebroid coproducts in $\Lc_-\rtimes\Ha$ and \tw{\Ha}{\Ru^+}{(\Lc_-\rtimes\Ha)}, when restricted to
$\Ha$ and respectively to $\Ha\tp \Ha$,
are obtained from the Hopf algebra coproducts
in $\Ha$ and \mbox{\tw{\Ha}{\Ru}{\Ha}}
by projecting the tensor products over $k$ to those over $\Lc$.
Recall that the coproduct of $\Ha$ defines a Hopf algebra map
from $\Ha$ to \mbox{\tw{\Ha}{\Ru}{\Ha}}.
It follows from here that $\eta$ respects
the bialgebroid coproducts when restricted to $\Ha$,
since $\eta$ is an $\Lc$-bimodule map.

Finally, it is obvious that $\eta$ respects the counits.
\end{proof}
\begin{remark}
\label{comp_twist}
If $\F$ is a dynamical twisting cocycle as in Proposition \ref{dc-tw},
then the twisted quantum groupoid  $\widetilde{\U\tp \Ha_\Lc}$
can be considered as a result of two consecutive twists
$$
\U\tp \Ha_\Lc\stackrel{\Psi_\Ru}{\longrightarrow}
\mbox{ \tw{\U}{\Ru}{\Ha_\Lc}}
\stackrel{\Psi_{\Ru^{-1}\F\Ru}}{\longrightarrow}
 \widetilde{\U\tp \Ha_\Lc}
$$
\end{remark}
\begin{corollary}
\label{twTPUDH}
Suppose that $\Ha$ is an arbitrary Hopf algebra and $\Lc$ is an $\Ha$-base algebra.
Let $\Ha$ be a Hopf subalgebra in a Hopf algebra $\U$.
Then there exists a bialgebroid homomorphism from $\Lc\rtimes \Ha$ to $\mbox{\tw{\U}{\Theta}{\D\Ha_\Lc}}$,
where $\Theta\in \Ha^*_{op}\tp \Ha$ is the standard quasitriangular structure on $\D\Ha$.
\end{corollary}
\begin{proof}
Recall from Remark \ref{rm+-} that $\Ha=(\D\Ha)^{\sst\{-\}}$ with respect to the standard
quasitriangular structure $\Theta\in \Ha^*_{op}\tp \Ha\subset (\D\Ha)^{\tp2}$.
Applying Proposition \ref{l_hom} to this case, we obtain the sequence of bialgebroid homomorphisms
$$
\Lc\rtimes\Ha  \to \Lc_-\rtimes\D\Ha\to
\mbox{\tw{\D\Ha}{\Theta}{(\Lc_-\rtimes\D\Ha)}}
\to\mbox{ \tw{\D\Ha}{\Theta}{\D\Ha_\Lc}},
$$
where the left arrow is embedding and the right one is projection along the ideal
$\D\Ha\tp J_\phi$.
The middle arrow is the map (\ref{hom}) where $\Ha$ is replaced by $\D\Ha$
and $\Ru$ by $\Theta$. This map is constructed out of the coaction
$\delta\colon\Lc\to \Ha\tp \Lc$ and the coproduct of $\D\Ha$. It remains to notice that,
as a coalgebra, the double $\D\Ha$ is a trivial tensor product of coalgebras
$\Ha$ and $\Ha^*_{op}$, whence the composite map takes the values
in $\mbox{\tw{\Ha}{\Theta}{\D\Ha_\Lc}}\subset \mbox{\tw{\U}{\Theta}{\D\Ha_\Lc}}$.
\end{proof}
\begin{remark}
\label{twTPop}
Replacing $\Ha$ by $\Ha_{op}$,  $\Lc$ by $\Lc_{op}$, and
$\Theta$ by $\bar \Theta=\Theta^{-1}$ in Corollary \ref{twTPUDH}, one
 can construct the twisted opposite bialgebroids
\tw{\U}{\bar \Theta}{(\Lc\rtimes\Ha)_{op}} and
\tw{\U}{\bar \Theta}{(\Ha_\Lc)_{op} },
cf. Definitions \ref{smash_op} and \ref{op_QGd}.
\end{remark}
\section{Dynamical categories and representations of bialgebroids}
\label{secFDCR}
In this section we establish relations between dynamical categories
from Definition \ref{defDC} and representations of bialgebroids
over quasi-commutative base. In this section $\Ha$ is
an arbitrary Hopf algebra and $\Lc$ is an $\Ha$-base algebra.

\subsection{Category $\Mod \Lc\rtimes \Ha$}
A central role in our further consideration belongs to Lemma \ref{U-L} below.
Let $X$ be a left $\Ha$-module.
Denote by  $\tilde X_\Lc$
the $\Lc$-bimodule $X\tp \Lc$
with respect to the following left and right action:
\be
\label{bimodule}
\la\llcorner(x\tp \mu)= \Theta_2\tr x\tp (\Theta_1\tr \la)\mu
,\quad
(x\tp \mu)\lrcorner\la= x\tp \mu\la
\ee
where $x\tp \mu\in \tilde X_\Lc$ and   $\la \in \Lc$.
For an $\Ha$-equivariant map $\psi\colon X\to Y\tp \Lc$ let $\tilde \psi_\Lc$ denote
the composition map
$$
X\tp \Lc\stackrel{\psi\tp \id_\Lc}{\longrightarrow}Y\tp \Lc\tp \Lc\stackrel{\id_Y\tp
\mathrm{m}_\Lc}{\longrightarrow} Y\tp \Lc,
$$
where $\mathrm{m}_\Lc$ is the multiplication in $\Lc$.
It is a morphism of $\Lc$-bimodules, due to the $\Ha$-invariance of $\psi$ and the
quasi-commutativity of $\Lc$.
\begin{lemma}
\label{U-L}
The correspondence $X\to \tilde X_\Lc$, $\psi\to \tilde \psi_\Lc$ commutes with
taking tensor products and defines a strong monoidal functor, $\Xi$,
from $\overline{\Mc}_{\Ha;\Lc}$ to the category  $\B i(\Lc)$ of $\Lc$-bimodules.
\end{lemma}
\begin{proof}
Straightforward.
\end{proof}
Denote by $\Modd\Lc\rtimes \Ha$ the full subcategory of $\Lc\rtimes \Ha$-modules
of the form $X\tp \Lc$, where $X$ is an $\Ha$-module.
\begin{thm}
\label{d-b}
The functor $\Xi$ establishes an isomorphism from the dynamical
category $\overline{\Mc}_{\Ha;\Lc}$ to $\Modd\Lc\rtimes \Ha$.
\end{thm}
\begin{proof}
Let $X\in \Ob\: \overline{\Mc}_{\Ha,\Lc}$ be an $\Ha$-module,
and $\tilde X_\Lc=\Xi(X)$ its image in $\B i(\Lc)$.
Consider $\tilde X_\Lc$ as an $\Ha$-module being the tensor product of $\Ha$-modules
$X$ and $\Lc$. One can check that this action together
with the left action of $\Lc$ gives rise to an action of $\Lc\rtimes \Ha$.
Further,  the tensor product of two $\Lc\rtimes\Ha$-modules $\tilde X_\Lc$ and $\tilde Y_\Lc$ is
$\Xi(X\tp Y)$, due to Lemma \ref{U-L}.

For  any morphism
$\psi\in \Hom_{\overline{\Mc}_{\Ha,\Lc}}(X,Y)$
the map $\tilde \psi_\Lc\colon \tilde X_\Lc\to \tilde Y_\Lc$ commutes with the action of  $\Lc\rtimes\Ha$.
Conversely, let $\phi\colon \tilde X_\Lc\to \tilde Y_\Lc$ be an  $\Lc\rtimes\Ha$-intertwiner.
Then $\phi$ is an $\Lc$-bimodule map and must have the form $\phi(x\tp \mu)=\phi(x\tp 1)\mu$.
The map $ \psi\colon X\to  \tilde X_\Lc$, $\psi(x):= \phi(x\tp 1)$,
is $\Ha$-equivariant, therefore it is a morphism in $\overline{\Mc}_{\Ha,\Lc}$ and
$\phi = \tilde \psi_\Lc$. Thus we have proved that the image of $\Xi$
is a full subcategory in $\Mod\Lc\rtimes \Ha$.
\end{proof}

Now suppose that $\Ha$ is a quasitriangular Hopf algebra
with the R-matrix $\Ru$ and $\Lc$ is $\Ha$-commutative.
Denote by $\Mc'_\Ha$ the full subcategory in $\Mc_\Ha$ consisting of such $\Ha$-modules $X$
that
\be
\Ru^+_{1}\tr x\tp \Ru^+_{2}\tr\la=\Ru^-_{1}\tr x\tp \Ru^-_{2}\tr\la
\label{cat'}
\ee
for all $x\tp \la\in X\tp \Lc$. Let $\overline{\Mc}'_{\Ha,\Lc}$ denote the dynamization
of $\Mc'_\Ha$,
i.e. the full subcategory in $\overline{\Mc}_{\Ha,\Lc}$ whose objects belong to $\Mc'_\Ha$.
Denote by $\Modd \Ha_\Lc$ the full subcategory of $\Ha_\Lc$-modules
of the form $X\tp \Lc$, where $X$ is an $\Ha$-module satisfying the condition (\ref{cat'}).
\begin{propn}
The category
$\overline{\Mc}'_{\Ha,\Lc}$ is a braided monoidal category. It
is naturally isomorphic to the category $\Modd \Ha_\Lc$,
which itself is a full subcategory in  $\Mod \Lc\rtimes\Ha$.
\end{propn}
\begin{proof}
The bialgebroid $\Ha_\Lc$ is the quotient of the bialgebroid
$\Lc\rtimes \Ha$ by the relations
$\Theta^+_2\tr \la\otimes \Theta^+_1=\Theta^-_2\tr \la\otimes \Theta^-_1$,
for all $\la \in \Lc$.
Therefore $\Mod \Ha_\Lc$ consists of those $\Lc\rtimes\Ha$-modules
whose annihilator contains this
ideal, thus $\Mod \Ha_\Lc$ is a  full subcategory in $\Mod \Lc\rtimes\Ha$.
An $\Lc\rtimes\Ha$-module $\Xi(X)$, where $X\in \Ob\: \overline{\Mc}_{\Ha,\Lc}$,
belongs to $\Mod \Ha_\Lc$ if and only if  $X\in \Ob\: \overline{\Mc}'_{\Ha,\Lc}$.
Applying Theorem \ref{d-b}, we conclude that restriction of the functor $\Xi$ to $\overline{\Mc}'_{\Ha,\Lc}$
gives an isomorphism of $\overline{\Mc}'_{\Ha,\Lc}$ with a full subcategory in $\Mod \Ha_\Lc$.
Since $\Mod \Ha_\Lc$ is braided, the category $\overline{\Mc}'_{\Ha,\Lc}$
is braided as well.
\end{proof}
\subsection{Category $\Mod\widetilde{\U\tp \D\Ha_\Lc}$}
In this subsection we assume that $\Ha$ is a Hopf subalgebra (not necessarily quasitriangular)
of a quasitriangular  Hopf algebra $\U$.
The category $\Mc_\U$ of $\U$-modules is viewed as a natural subcategory in $\Mc_\Ha$.

Any  module over the tensor product bialgebroid $\U\tp \D\Ha_\Lc$ can be represented as the tensor product
 $V\tp A$, where $V$ is an $\U$-module and $A$ is an $\D\Ha_\Lc$-module.
The induced $\Lc$-bimodule structure on $V\tp A$ coincides with the standard one:
$$
\la\llcorner(v\tp a):= v\tp s(\la)a,
\quad
(v\tp a)\lrcorner\la:= v\tp t(\la)a,
$$
for $v\tp a\in V\tp A$ and $\la\in \Lc$.

Let $\Lc$ be an $\Ha$-base algebra $\Lc$ and $\F$ a dynamical twist.
Consider the twisted bialgebroid $\widetilde{\U\tp \D\Ha_\Lc}$ built by
means of twist $\Psi=\Psi_{\F\Theta}$ from Subsection \ref{subsecTDC}.
\begin{propn}
Objects $V\tp \Lc$, where $V$ is a $\U$-module,
form a full monoidal subcategory, $\Modd\widetilde{\U\tp \D\Ha_\Lc}$,
in $\Mod\widetilde{\U\tp \D\Ha_\Lc}$. It is isomorphic
to $\Mc_\U$ if and only if $\Lc$ is quasi-transitive.
In the particular case of the unit $\F$, the
isomorphism from $\Mc$ to $\Modd \mbox{\tw{\U}{\Theta}{\D\Ha_\Lc}}$
is enclosed in the commutative diagram
\be
\begin{array}{ccc}
\Mc_\U & \longrightarrow & \overline{\Mc}_{\Ha,\Lc}
\\[6pt]
\wr \parallel &&\wr \parallel \Xi
\\[6pt]
\Modd \mbox{\tw{\U}{\Theta}{\D\Ha_\Lc}} & \longrightarrow &
\Modd \Lc\rtimes \Ha
\end{array},
\ee
where the bottom line is induced by the bialgebroid homomorphism
$\mbox{\tw{\U}{\Theta}{\D\Ha_\Lc}}\longleftarrow \Lc\rtimes \Ha$.
\end{propn}
\begin{proof}
The $\U\tp \D\Ha_\Lc$-modules of the form $V\tp \Lc$, where
$V$ is a $\U$ module, are closed under the tensor product
induced by the twist. Thus they form a full monoidal subcategory
in $\Mod \U\tp \D\Ha_\Lc$.
Obviously,
$\Hom_{\U\tp \D\Ha_\Lc}(V\tp \Lc, W\tp\Lc)\simeq \Hom_{\U}(V,W)\tp\End_{\D\Ha}(\Lc)$.
There is a natural bijection between $\End_{\D\Ha}(\Lc)$ and $\Lc^{\D\Ha}$,
following from Lemma \ref{lemcenter}.
Thus the category $\Mc_\U$ is isomorphic to $\Modd\widetilde{\U\tp \D\Ha_\Lc}$, provided
$\Lc$ is quasi-transitive.
In the particular case of the unit $\F$, the isomorphism
from $\Mc_\U$ to $\Modd \mbox{\tw{\U}{\Theta}{\D\Ha_\Lc}}$
is given by the functor $\Xi$.
\end{proof}

Summarizing, we present a diagram of most important bialgebroids
and their interrelations:
$$
\begin{array}{cccccc}
\U\tp\D\Ha_\Lc
\\[6pt]
\hspace{5pt}\wr\hspace{5pt}\Psi_{\Theta}
\\[6pt]
\mbox{\tw{\U}{\Theta}{\D\Ha_\Lc}}&\leftarrow  \Lc\rtimes\Ha\rightarrow \D\Ha_\Lc\;,
\\[6pt]
\hspace{23pt}\wr\hspace{5pt}\Psi_{\Theta^{-1}\F\Theta}
\\[6pt]
\widetilde{\U\tp\D\Ha_\Lc}
\end{array}
$$
The columns represent the bialgebroid  twist, cf. Remark \ref{comp_twist}.
The horizontal arrows are bialgebroid homomorphisms.

Assuming the $\Ha$-base algebra $\Lc$ to be quasi-transitive
and passing to the $\Modd$-modules, we obtain the following commutative diagram
displaying the interrelations between the categories:
$$
\begin{array}{ccccccccc}
\Modd (\U\tp\D\Ha_\Lc)\simeq \Mc_\U&\longrightarrow&\overline{\Mc}_{\Ha,\Lc}&
\leftarrow&\overline{\Mc}_{\Ha,\Lc}'
\\[6pt]
\wr\parallel\Psi_{\Theta}\hspace{15pt}&&\wr\parallel\Xi\hspace{10pt}&&\wr\parallel\Xi\hspace{10pt}
\\[6pt]
\Modd \mbox{\tw{\U}{\Theta}{\D\Ha_\Lc}}&
\longrightarrow & \Modd\Lc\rtimes\Ha&\leftarrow &\Modd\D\Ha_\Lc
\\[6pt]
\hspace{10pt}\wr\parallel\Psi_{\Theta^{-1}\F\Theta}
\\[6pt]
\Modd \widetilde{\U\tp\D\Ha_\Lc}.
\end{array}
$$
\section{Dual quantum groupoids (dynamical FRT algebras)}
\label{secDQG}
In this section we present an example of a module algebra over the bialgebroid $\widetilde{\U\tp\D\Ha_\Lc} $ constructed
in Subsection \ref{subsecTDC}.
It turns out to be  a bialgebroid and  may be thought of as an analog of the dual
Hopf algebra. This fact is not occasional and will be addressed in a separate
publication.
\subsection{Dynamical associative algebras}

Let $\Ha$ be a Hopf algebra and $\Lc$ its base algebra. Recall that $\overline{\Mc}_{\Ha,\Lc}$
denotes the dynamical extension  over $\Lc$ of the category $\Mc_\Ha$ of left $\Ha$-modules.
\begin{definition}[\cite{DM1}]
\label{defDA}
\select{Dynamical associative  algebra} (or simply dynamical algebra)
is an algebra in the monoidal category $\overline{\Mc}_{\Ha,\Lc}$.
\end{definition}
A dynamical algebra $\A$ is an $\Ha$-module equipped
with an $\Ha$-equivariant bilinear map $\divideontimes\colon\A\tp
\A\to \A\tp \Lc$ such that the following diagram is commutative:
\be
&
\begin{diagram}
   \dgARROWLENGTH=0.5\dgARROWLENGTH
\node{\A \tp \Lc\tp \A}\arrow{e,t}{\id \tp\tau_\A}\node{\A \tp \A \tp \Lc}
\arrow{e,t}{\divideontimes\tp \id}\node{\A \tp \Lc \tp \Lc}\arrow{e,t}{\id\tp\mathrm{m}_\Lc}
\node{\A \tp \Lc} \arrow{s,r,!}{\parallel}
\\
\node{\A\tp\A\tp \A }\arrow{e,t}{\id\tp \divideontimes }\arrow{n,l}{\divideontimes\tp \id}
\node{\A \tp \A \tp \Lc}\arrow{e,t}{\divideontimes\tp \id}\node{\A \tp \Lc \tp \Lc}\arrow{e,t}{\id\tp\mathrm{m}_\Lc}\node{\A \tp \Lc}
\end{diagram}
\label{shifted_assoc0}
\ee
Here $\mathrm{m}_\Lc$ stands for the multiplication in $\Lc$ and
$\tau_\A$ denotes the permutation $\Lc\tp \A\to  \A\tp \Lc$, which
is expressed through the coaction $\delta(\la)=\la^{(1)}\tp \la^{[2]}$
or, equivalently, through the canonical R-matrix of the double $\D\Ha$ by
$\la\tp a\mapsto \la^{(1)}\tr a\tp \la^{[2]}=\Theta_2\tr a\tp \Theta_1\tr \la$.

If the operation $\divideontimes$ takes  values in $\A\tp 1\subset \A\tp\Lc$,
the condition (\ref{shifted_assoc0}) reduces to the ordinary associativity.
For example, suppose that  $\A$ is a module algebra over a Hopf algebra $\U$ containing $\Ha$,
then it is a dynamical algebra. Let $\F=\F_1\tp\F_2\tp \F_3\in \U\tp\U\tp \Lc$ be
a dynamical cocycle.
Then the map $\A\tp\A\to \A\tp \Lc$,
\be
\label{twised_alg}
a\divideontimes b:= (\F_1\tr a)(\F_2\tr b) \tp \F_3,
\ee
defines a new structure of dynamical algebra on $\A$, \cite{DM1}.
\begin{propn}
\label{da-abi}
A left $\Ha$-module $\A$ is a dynamical associative algebra if and only if
$\A\tp \Lc$ equipped with the $\Lc$-bimodule structure (\ref{bimodule}) is an algebra
in the category of $\Lc$-bimodules.
\end{propn}
\begin{proof}
This follows from the isomorphism of categories $\overline{\Mc}_{\Ha,\Lc}\simeq \Modd \Lc\rtimes \Ha$.
\end{proof}
\begin{remark}
\label{alg-alg}
Remark that an $\Lc$-bimodule and associative algebra
$\B$ is an algebra in the category of $\Lc$-bimodules
if and only if the multiplication and $\Lc$-actions are compatible:
 $(a\lrcorner \la)b=a(\la\llcorner b)$ for all $a,b\in \B$ and $\la\in \Lc$.
\end{remark}
We will denote the associative algebra from proposition \ref{da-abi} by  $\A\ltimes \Lc$

\subsection{Bialgebroid $\U^*_\F\ltimes(\Lc\tp \Lc_{op})$}
\label{secDual}
In this subsection we construct a dynamical dual to a Hopf algebra $\U$.
The dynamical dual can be regarded as a dynamical analog of the FRT algebra
if $\U$ is quasitriangular.

We consider the dual Hopf algebra $\U^*$ as a $\U\tp \U_{op}$-module,
with respect to the coregular actions
\be
\label{D-inv}
x\tr u:= u^{(1)}\langle x,u^{(2)}\rangle, \quad y\btr u:= \langle y,u^{(1)}\rangle u^{(2)},
\ee
where $x\in \U$, $y\in \U_{op}$, $u\in \U^*$. By restriction,  $\U^*$ is also
an $\Ha\tp \Ha_{op}$-module algebra.
Recall from Lemma \ref{op} that $\Lc_{op}$ is an $\Ha_{op}$-base algebra.
By this reason, we can consider $\U^*$ as a dynamical associative algebra over the $\Ha\tp \Ha_{op}$-base algebra
$\Lc\tp \Lc_{op}$.
Applying Proposition \ref{da-abi}, we construct an algebra, $\U^*\ltimes(\Lc\tp \Lc_{op})$,
in the category of $\Lc\tp \Lc_{op}$-bimodules. Moreover, it is an algebra in
the category of modules over
$\bigl(\mbox{\tw{\U}{\Theta}{\Ha_\Lc}}\bigr)\tp\bigl(\mbox{\tw{\U_{op}}{\bar\Theta}{(\Ha_\Lc)_{op}}}\bigr)$,
which is a tensor product of $\Lc$- and $\Lc_{op}$-bialgebroids,
cf. Example \ref{H-B} and Remark \ref{twTPop}.
By Remark \ref{alg-alg}, $\U^*\ltimes(\Lc\tp \Lc_{op})$
is an associative algebra with the multiplication
$$
(u\tp \la\tp \mu)(v\tp \al\tp \bt):=u \bigl(\la^{(1)}\tr \mu^{(1)}\btr v\bigr)\tp  \la^{[2]}\al\tp \mu^{[2]}\cdot\bt
$$
for $u,v\in \U^*$, $\la,\al\in \Lc$ and $\mu,\bt\in \Lc_{op}$.

The following two propositions can be checked by direct but tedious calculations.
\begin{propn}
The algebra $\U^*\ltimes(\Lc\tp \Lc_{op})$ is a right $\Lc$-bialgebroid (cf. Remark \ref{right}) over the base $\Lc$
with the target map $t\colon \la\mapsto 1\tp 1\tp \iota(\la)$, the source map $s\colon \la\mapsto 1\tp \la\tp 1$,
the coproduct $\Delta\bigl(u\tp \la\tp \mu\bigr):=
\bigl(u^{(1)}\tp 1\tp \mu\bigr)\tp_{\Lc}(u^{(2)}\tp \la\tp 1)$,
and the counit
$\ve\bigl(u\tp \la\tp \mu\bigr):=
\ve_{\U^*}(u)\: \iota^{-1}(\mu)\la
$.
\end{propn}

Suppose that $\F\in \U\tp \U\tp \Lc $ is a universal dynamical cocycle
over $\Lc$.
Then one can check that
$\bar\F:=\F^{-1}\in \U_{op}\tp \U_{op}\tp \Lc_{op}$ is a universal dynamical cocycle
over $\Lc_{op}$, which is an $\Ha_{op}$-base algebra.
Thus $\F\tp \bar\F$ is a universal twist in the dynamical extension of the category
of $\U$-bimodules over the base algebra $\Lc\tp \Lc_{op}$.
Let $\U^*_\F$ denote  the dynamical associative algebra over
the base $\Lc\tp \Lc_{op}$ obtained from
$\U^*$ by the twist $\F\tp \bar\F$, see (\ref{twised_alg}).
\begin{propn}
The algebra  $\U^*_\F\ltimes(\Lc\tp \Lc_{op})$ is
a right $\Lc$-bialgebroid with the same source and target maps, coproduct, and counit
as  $\U^*\ltimes(\Lc\tp \Lc_{op})$.
It is an algebra in the category
of modules over the $\Lc\tp \Lc_{op}$-bialgebroid
$\widetilde{\U\tp \Ha_\Lc}\tp \widetilde{\U_{op}\tp (\Ha_\Lc)_{op}}$.
\end{propn}
\begin{remark}
The dynamical algebra  $\U^*_\F$ is obtained from  $\U^*$ by the dynamical
twist $\F$ applied from the two sides.
Applied to only one side, the dynamical twist gives a dynamical algebra
in the category $\overline{\Mc}_{\Ha,\Lc}$, which participates in the
equivariant star product quantization of vector bundles on coadjoint orbits
of reductive Lie groups, \cite{DM1}.
\end{remark}

\section{On the quasi-classical limit and dynamical r-matrix}
\label{secCDRM}
In this section we consider Lie bialgebroids that are relevant
to quantum groupoids studied in this paper.
For an exposition of the theory of Lie algebroids and Lie bialgebroids,
the reader is referred to \cite{K-Schw} and \cite{MXu}.
\subsection{Lie bialgebroids}
Let us recall that a Lie algebroid $\Bg_0$ over a commutative algebra  (sheaf) $\Lc_0$ is
an $\Lc_0$-module equipped with a structure of Lie algebra together
with a Lie algebra homomorphism (anchor) $\Bg_0\to Der(\Lc_0)$ such that
$[\xi,f\eta]=f[\xi,\eta]+(\xi\tr f)\;\eta$ for all $\xi,\eta \in \Bg_0$,
$f\in \Lc_0$.
\begin{example}
\label{trivial}
Let $\h$ be a Lie algebra acting on $\Lc_0$. Consider the trivial bundle
$\Lc_0\tp \h$ equipped with the following Lie algebra structure on section:
$$[f\tp \xi,g\tp \eta]:=fg\tp [\xi,\eta]+ f(\xi\tr g) \tp \eta - g(\eta\tr f) \tp \xi$$
for $f,g\in \Lc_0$ and $\xi,\eta\in \h$. The anchor map is determined by the
action of $\h$ on $\Lc_0$.
We denote this Lie bialgebroid by $\Lc_0\rtimes \h$.
\end{example}
The Lie bracket on $\Bg_0$ can be extended as the Schouten bracket $[.,.]$
to the exterior algebra $\wedge^\bullet \Bg_0$ making it a Gerstenhaber algebra,
see \cite{K-Schw}. A Lie algebroid structure defines a "de Rham" differential $d$
of degree 1 with zero square on the graded exterior  algebra
$\wedge^\bullet \Bg^*_0$ of the dual vector bundle $\Bg^*$.

With every Lie algebroid $\Bg_0$ one can associate a universal enveloping $\Lc_0$-bialgebroid, $\U(\Bg_0)$, see \cite{Xu1}.
For the Lie algebroid from Example \ref{trivial}, it coincides with $\Lc_0\rtimes \U(\h)$,
where  $\U(\h)$ is the universal enveloping algebra of $\h$ (note that
$\Lc_0$ is a base algebra for $\U(\h)$), in the sense of Definition \ref{defBA}.

Infinitesimal theory of quantized universal enveloping algebras leads to the notion of Lie bialgebras.
Analogously, the problem of quantization of $\U(\Bg_0)$ in the class of bialgebroids gives rise to
the notion of Lie bialgebroid, \cite{MXu}. By definition, a Lie bialgebroid
is a pair $(\Bg_0,\Bg_0^*)$ of two Lie algebroids in duality
satisfying the compatibility condition
$$
d_*[\xi,\eta]=[d_*\xi,\eta]+[\xi,d_*\eta].
$$
Below we give examples of Lie bialgebroids that are relevant to our study.

Let $(\h,\h^*)$ be a Lie bialgebra and $\D\h:=\h\bowtie \h_{op}$ its double Lie (bi)algebra.
\begin{definition}[\cite{DM1}]
A Poisson $\h$-base algebra $\Lc_0$ is a commutative
algebra equipped with a left $\D \h$-action such that the canonical symmetric invariant tensor
$\theta\in \D\h^{\tp 2}$ vanishes on $\Lc_0$.
\end{definition}
It follows from the definition that the classical r-matrix of $\D \h$ induces a Poisson
bivector field on $\Lc_0$. This Poisson structure can be quantized to a $\U_q(\h)$-base algebra $\Lc$ for $\U_q(\h)$
being the quantization of $\U(\h)$ along $\h^*$, \cite{DM1}.

One can check the following
\begin{propn}
\label{Lie-b1}
Let $\h$ be a Lie algebra and $\h^*$ a Lie algebra structure on the dual space. The Lie algebroids $\Lc_0\rtimes \h$ and $\Lc_0\rtimes \h^*$
form a Lie bialgebroid {\em iff} $(\h,\h^*)$ is a Lie bialgebra and $\Lc_0$ is a
Poisson $\h$-base algebra.
The differential $d_*$ is given by the Lie cobracket $\nu$ on $\h$ considered
as a constant section of $\wedge^2(\Lc_0\rtimes \h)$.
\end{propn}

A Lie bialgebroid $\Bg_0$ is called coboundary if the differential $d_*$ is
generated by an element $\Lambda\in \wedge^2 \Bg_0$,
namely, has the form $d_*(\zeta):=[\Lambda,\zeta]$, where $\Lambda$
satisfies the condition $[[\Lambda,\Lambda],\Lambda]=0$.
The Lie bialgebroid  $\Lc_0\rtimes \h$ is not coboundary, in general,
even in the case of coboundary Lie bialgebra $\h$. However, if $\h$ is quasitriangular,
then there is an ideal in $\Lc_0\rtimes \h$ lying in the kernel of the
anchor map. The quotient of $\Lc_0\rtimes \h$ by that ideal is a coboundary
Lie bialgebroid.
We will demonstrate this on the example of $\Lc_0\rtimes \D\h$
(note that a Poisson $\h$-base algebra is that for $\D\h$ as well).
\begin{example}
\label{Lie-b2}
Suppose that $\Lc_0$ is a Poisson $\h$-base algebra.
Let $\{h_i\}$ be a base in $\h$ and $\{\eta^i\}$ its dual in $\h^*_{op}$.
Let $\theta:=\frac{1}{2}\sum_i (h_i\tp\eta^i+\eta^i\tp h_i)$ denote
the canonical symmetric invariant element of the double
$\D\h$.
Consider in $\Lc_0\rtimes \D\h$ an $\Lc_0$-submodule generated by
sections of the form $\theta_1\tr f\tp \theta_2$ for all $f\in \Lc_0$.
It forms
an ideal $J_0$ in the Lie algebra $\Lc_0\rtimes \D\h$, and this ideal
is  $\D\h$-invariant.
The quotient of $\Lc_0\rtimes \D\h$ by $J_0$ is a coboundary Lie algebroid, $\D\h_{\Lc_0}$.
Its dual Lie bialgebroid is the  annihilator of $J_0$  in $\Lc_0\rtimes \D^*\h$.

Suppose now that $\Lc_0$ is a function algebra on a Poisson $\h$-base manifold $L$, which
is assumed $\D\h$-homogeneous.
Then $J_0$ can be considered as the space of sections of an $\h$-vector bundle over $L$.
Let us fix an origin in $L$ and let $\mathfrak{k}\in \h$ be the Lie
algebra of its stabilizer.  Denote by $\mathfrak{k}_0$ the ideal in $\mathfrak{k}$ that is the
kernel of the canonical invariant inner product in $\D\h$ restricted to $\mathfrak{k}$.
Then the fiber of $J_0$ is isomorphic to $\mathfrak{k}_0$.
\end{example}

Given a Lie bialgebroid over $\Lc_0$, the latter is equipped with a Poisson structure
$\{f,g\}:=(df,d_*g)$. Quantization of a Lie bialgebroid $\Bg_0$ over $\Lc_0$
means quantization, $\Lc_\hbar$, of $\Lc_0$ and construction of an $\Lc_\hbar$-bialgebroid
whose a) classical limit is the universal enveloping
$\Lc_0$-bialgebroid $\U(\Bg_0)$ and b) the infinitesimal deformation is determined
by the Lie bialgebroid $\Bg_0$.
Conversely, an $\Lc_\hbar$-bialgebroid
$\Bg_\hbar$ over $k=\C[[\hbar]]$ such that $\Lc_0=\Lc_\hbar\mod \hbar$ and $\Bg_\hbar=\U(\Bg_0)\mod \hbar$
gives rise to a structure of Lie bialgebroid over $\Lc_0$ in the quasi-classical limit, \cite{MXu}.

Let $\h$ be a Lie bialgebra, $\Lc_0$ a Poisson base algebra over $\h$,
$\U_\hbar(\h)$ and $\U_\hbar(\D\h)$ the corresponding quantizations of the universal enveloping algebras,
and $\Lc$ is the $\U_\hbar(\h)$-base algebra that is a quantization of $\Lc_0$.
Then the  $\Lc$-bialgebroid $\Lc\rtimes\U_\hbar(\h)$ is a quantization of the Lie bialgebroid
$\Lc_0\rtimes \h$ from Proposition \ref{Lie-b1}.
The ideal $J_0$ from Example \ref{Lie-b2} is a classical limit of the biideal $J_\phi$
from Proposition \ref{q-triang}, where the role of the double belongs to
the arbitrary quasitriangular Hopf algebra.
The quantum groupoid  $\U_\hbar(\D\h)_\Lc$
is a quantization of  the coboundary Lie bialgebroid
and $\D\h_{\Lc_0}$ from Example \ref{Lie-b2}.

In the next subsection we describe more complicated coboundary Lie bialgebroids,
which are related to dynamical r-matrices. The corresponding theory
for commutative base was developed by Xu. Here we consider an arbitrary
Lie bialgebra $\h$ and its base manifold $L$. It turns out that
the dynamical r-matrices are in one-to-one correspondence with
a class of Lie bialgebroid structures on certain Lie algebroids.
\subsection{Classical dynamical r-matrix and Lie bialgebroids}
Let $\h$ be a Lie bialgebra, with the Lie cobracket $\nu\colon \h\to \h\wedge \h$.
Let $\Lc_0$ be a Poisson $\h$-base algebra. The reader may think of $\Lc_0$ as
a function algebra on a base manifold.

Suppose that $\h$ is a subalgebra in
a Lie algebra $\g$.
We emphasize that we do not assume {\em any} Lie bialgebra structure on $\g$.
The algebra $\Lc_0$ is equipped with an $(\g\oplus \D\h)$-action
assuming it trivial when restricted to $\g$.
\begin{definition}[DM1]
\label{defDrm}
An element  $
r(\la) \in \Lc_0\tp\wedge^2\g$ is called a \select{classical dynamical r-matrix}
over the Poisson base algebra $\Lc_0$ with values in $\wedge^2 \g$ if
\begin{enumerate}
\item for any $h\in\h$
\be \label{inv}
 h\tr   r(\la) + [h\tp 1+ 1\tp h,  r(\la)] = \nu(h),
\ee
\item
$  r$ satisfies  the  equation
\be
\label{cDYBE} \sum_{i}
\mathrm{Alt}\bigl(h_i\tp  {\eta}^i \tr  r(\la)\bigr)-\mathrm{CYB}\bigl(r(\la)\bigr)
 = \varphi(\la) \in \Lc^{\D\h}_0\tp(\wedge^3\g)^{\g}=(\Lc_0\tp\wedge^3\g)^{\g\oplus \D\h},
\ee
where $\mathrm{CYB}(\zeta):=[\zeta_{12},\zeta_{13}]+[\zeta_{13},\zeta_{23}]+[\zeta_{12},\zeta_{23}]$,
$\zeta\in \Lc_0\tp \wedge^2 \g$, is the Yang-Baxter operator,
$\mathrm{Alt}(\xi_1\tp\xi_2\tp\xi_3):=\xi_1\tp\xi_2\tp\xi_3-\xi_2\tp\xi_1\tp\xi_3+\xi_2\tp\xi_3\tp\xi_1$,
$\xi_i\in \Lc_0\tp \g$,
and $\varphi(\la)$ is some invariant element.
\end{enumerate}
\end{definition}

Consider the trivial Lie bialgebroid $\Lc_0\tp\g$ with the zero anchor map. It is
just a Lie algebra over $\Lc_0$.
Denote by $(\g\oplus \D\h)_{\Lc_0}$ the direct sum Lie algebroid $(\Lc_0\tp\g)\oplus(\D\h_{\Lc_0})$.
Let $\{h_i\}$ be a base in $\h$ and $\{\eta^i\}$ its dual in $\h^*_{op}$.
Consider the sum $\Lambda_0:=\varpi+2\varpi'\in \wedge^2(\g\oplus \D\h)$, where
$\varpi:=\sum_i \eta^i\wedge h_i=\frac{1}{2}\sum_i (\eta^i\tp h_i-h_i\tp \eta^i)\in \wedge^2 \D\h$ denotes the universal r-matrix
of the double, and $\varpi'\in \g\wedge \h^*_{op}$
is obtained from $\varpi$ via the embedding $\h\to \g$.
The element $\Lambda_0$ can be thought of as a constant section
of the exterior square of the trivial vector bundle $(\Lc_0\tp\g)\oplus(\Lc_0\rtimes\D\h)$.
Let us denote the projection of $\Lambda_0$ to
$\wedge^2 (\g\oplus \D\h)_{\Lc_0}$ by the same letter.

\begin{thm}
\label{r-L-b}
a) Let $r(\la)\in\Lc_0\tp \wedge^2\g$ be  a classical dynamical r-matrix,
Then the element $\Lambda\!:=  r(\la)+\Lambda_0 \in \wedge^2 (\g\oplus \D\h)_{\Lc_0}$,
generates a zero square differential on the graded Lie algebra $\wedge^\bullet(\g\oplus \D\h)_{\Lc_0}$ and therefore defines
a coboundary Lie bialgebroid on $(\g\oplus \D\h)_{\Lc_0}$.

b) Suppose that $\h^*_{op}$ acts effectively on $\Lc_0$ and
the element $\Lambda=  r(\la)+\Lambda_0 \in \wedge^2 (\g\oplus \D\h)_{\Lc_0}$ defines a coboundary
Lie bialgebroid on $(\g\oplus \D\h)_{\Lc_0}$. Then  $r(\la)\in\Lc_0\tp \wedge^2\g$
is a classical dynamical r-matrix.
\end{thm}
\begin{proof}
The element $\Lambda$ defines a coboundary Lie bialgebroid on $(\g\oplus \D\h)_{\Lc_0}$
if and only if
\be
[[\Lambda,\Lambda],f\tp\xi] =0\mod J_0, \quad \forall f\in \Lc_0,\> \forall\xi\in \g\oplus \D\h.
\ee
Explicitly,
the ideal $J_0$ is generated by the relations
\be
\label{J0}
\theta_1\tr f\tp \theta_2= \frac{1}{2}\sum_i  \bigl((\eta^i\tr f )\tp h_i+ (h_i\tr f )\tp \eta^i\bigr )=0
\ee
for all $f\in \Lc_0$, see Example \ref{Lie-b2}.

Let us denote by $\mathrm{Span}(. )$ the $\Lc_0$-module generated by a given set.
We will analyze the structure of the Schouten bracket $[\Lambda,\Lambda]\in \mathrm{Span}\wedge^{3}(\g\oplus \D\h)$.

1) The contribution of $[\Lambda,\Lambda]$ to
$\mathrm{Span} \bigl(\g\wedge \D\h\wedge \D\h \bigr)$ is proportional to
$$
-h_i\wedge[\eta^i ,\eta^j]\wedge h_j- h_i\wedge\eta^j\wedge [\eta^i ,h_j]+[h_i, h_j]\wedge \eta_i\wedge \eta_j ,
$$
where the summation over repeating indices
is understood. This term is identically zero, which
follows from definition of the double.

2) The contribution of  $[\Lambda,\Lambda]$ to
$\mathrm{Span} \bigl(\wedge^3\D\h\bigr)$ is equal to $[\varpi,\varpi]$ that is proportional to
$[\theta_{12},\theta_{23}]$.
The latter is a $\g\oplus \D\h$-invariant, and $[[\varpi,\varpi],f]$ belongs
to $J_0\wedge \D\h\wedge\D\h$ for all $f\in \Lc_0$, see relations (\ref{J0}). It follows
that the commutator of $[\varpi,\varpi]$ with all elements of $\wedge^\bullet (\g\oplus \D\h)_{\Lc_0}$
vanishes.

3) The contribution of  $[\Lambda,\Lambda]$ to $\mathrm{Span} \bigl(\wedge^3\g\bigr)$ is equal
to
$$
\frac{4}{3}\varphi(\la):=[r(\la),r(\la)]-4 \sum_i h_i\wedge \bigl(\eta^i\tr r(\la)\bigr).
$$
This definition of $\varphi(\la)$ is equivalent to (\ref{cDYBE}).
The element $\varphi$ commutes with all elements from $\Lc_0$.
Its commutator with $\g\oplus \D\h$ cannot belong to $J_0\wedge (\ldots)$ and hence vanishes
if and only if $\varphi$ is $(\g\oplus \D\h)$-invariant.

4) In fact, the contribution of  $[\Lambda,\Lambda]$ to $\mathrm{Span} \bigl(\D\h\wedge \g\wedge\g\bigr)$
lies, modulo $J_0$, in $\mathrm{Span} \bigl(\h^*_{op}\wedge \g\wedge\g\bigr)$.
Using the identity (\ref{J0}) we find  this contribution to be  proportional, modulo $J_0$,
to
\be
\label{auxi_term}
\eta^i\wedge [h_i,r(\la)]+
\eta^i\wedge (h_i\tr r(\la))+
[\eta^i,\eta^j]\wedge h_i\wedge h_j
\ee
(summation understood).
If $r(\la)$ is a classical dynamical r-matrix, this term is identically zero.
This follows from equation (\ref{inv}), since $\sum_i[\eta^i,\eta^j]\wedge  h_i\wedge h_j=-\sum_i\eta^i\wedge \nu(h_i)$
(recall
that the commutator is taken in the opposite Lie algebra $\h^*_{op}$.

Thus we have proven that $\Lambda$  defines a structure of coboundary bialgebroid
if $r(\la)$ is a dynamical r-matrix.

Conversely, let $\Lambda$ define a coboundary bialgebroid on $(\g\oplus \D\h)_{\Lc_0}$.
The steps 1)-3) of the above proof imply that $r(\la)$ satisfies the equation (\ref{cDYBE}).
Assume now that  $\h^*_{op}$ acts effectively on $\Lc_0$.
Taking  commutator of (\ref{auxi_term}) with an arbitrary element $f\in \Lc_0$,
we find that the expression
\be
(\eta^i\tr f)\tp  [h_i,r(\la)]+
(\eta^i \tr f)\tp (h_i\tr r(\la))+
([\eta^i,\eta^j]\tr f)\tp  h_i\wedge h_j
\ee
(summation understood) is equal to zero if and only if the equation (\ref{inv}) is satisfied for
all $h\in \h$.
This completes the proof.
\end{proof}
\begin{remark}
The element $\varphi(\la)$ in the right-hand side of (\ref{cDYBE}) is constant,
i.e. belongs to  $(\wedge^3\g)^\g$
in case $\Lc_0$ is quasi-transitive, i.e. $\Lc_0^{\D\h}=k$.
\end{remark}

Equation (\ref{cDYBE}) with the zero right-hand side
is called classical dynamical Yang-Baxter equation (CDYBE) over the Poisson base algebra $\Lc_0$.
When the invariant element $\varphi$ is non-zero, it may be called
modified CDYBE.

Suppose that the element $\varphi(\la)$ can be resolved by a symmetric element
$ \omega(\la)\in \Lc_0^{\D\h}\tp\g^{\tp2}$
in the sense of the equality $\varphi(\la)=-[\omega_{12}(\la),\omega_{23}(\la)]$. Then
the element $  r(\la)+\omega(\la)$ will satisfy equation (\ref{inv})
and equation (\ref{cDYBE}) with zero $\phi$, although it will not
be skew-symmetric. Conversely, if an element $  r(\la)\in \Lc_0\tp\g\tp \g$
with $\g$-invariant symmetric part
satisfies equations  (\ref{inv}) and (\ref{cDYBE}), then its skew-symmetric
part is a dynamical r-matrix in the sense of Definition \ref{defDrm}.

The classical dynamical r-matrices were conventionally defined
on a "flat" base manifold,
\cite{BDFh,F,EV1,ESch2,Sch}, namely the dual space $\h^*$ with Lie algebra structure.
This corresponds to  the zero right-hand side of equation (\ref{inv}).
A lot of progress in quantization of such r-matrices, including
the Alekseev-Meinrenken solution \cite{AM} and its generalizations, \cite{ESch2}, has been made
in recent the recent papers of
Enriquez and Etingof, \cite{EE1,EE2}.

Dynamical (non-skew) r-matrices over an arbitrary Lie bialgebra $\h$ and
an $\h$-base manifold were introduced in \cite{DM1}. The definition
of dynamical r-matrix given in \cite{DM1} was slightly less general then in
the present paper. Namely, $\g$ was assumed to be a Lie bialgebra
containing $\h$ as a sub-bialgebra.
An example of dynamical r-matrix on a group manifold was given in \cite{FhMrsh}.
The existence of such r-matrices for a wide class of Lie bialgebras
follows from the fusion procedure of \cite{DM1} adopted to the quantum
group case.

Thus there arises a problem of classification of dynamical r-matrices
on non-flat base manifolds and the problem of their quantization.
In view of Theorem \ref{r-L-b}, the second problem
is closely related to the problem of quantization of Lie bialgebroids
of a special class.

\end{document}